\documentclass[a4paper,11pt]{article}

\usepackage{pdfsync}
\usepackage{calrsfs} 

\usepackage{mathrsfs}
\usepackage{amsmath}
\usepackage[applemac]{inputenc}
\usepackage{amsfonts}
\usepackage{amssymb}
\usepackage{amsthm}
\usepackage{stmaryrd}
\usepackage{color}
\usepackage{mathabx}
\usepackage{dsfont}

\usepackage{hyperref}
\usepackage{mathtools}                                                      
\mathtoolsset{showonlyrefs=true}                                    

\def\cC{{\mathcal{C}}}

\renewcommand{\theequation}{\thesection.\arabic{equation}}

\let\counterwithin\relax
\usepackage{chngcntr}

\counterwithin{equation}{section}

\newcommand{\Q}{\mathbf{d}}

\newcommand{\Z}{\mathbb{Z}}

\newcommand{\N}{\mathbb{N}}

\newcommand{\R}{\mathbb{R}}

\renewcommand{\theequation}{\thesection.\arabic{equation}}

\def\1{\mbox{1\hspace{-0.25em}l}}

\def\DB{{\mathfrak D}}
\def\Rdd{\mathbb{R}^{N+1}}

\def\bx{{\mathbf{x}}}
\def\x{{\boldsymbol{x}}}
\def\p{{\mathbf{p}}}

\def\y{{\eta}}
\def\l{\lambda}
\def\<{\langle}
\def\>{\rangle}

\def\x{{\xi}}
\def\z{{\zeta}}

\numberwithin{equation}{section}
\def\theequation{\arabic{section}.\arabic{equation}}
\newcommand{\be}{\begin{eqnarray}}
\newcommand{\ee}{\end{eqnarray}}
\newcommand{\ce}{\begin{eqnarray*}}
\newcommand{\de}{\end{eqnarray*}}
\newtheorem{theorem}{Theorem}[section]
\newtheorem{lemma}[theorem]{Lemma}
\newtheorem{remark}[theorem]{Remark}
\newtheorem{definition}[theorem]{Definition}
\newtheorem{proposition}[theorem]{Proposition}

\newtheorem{corollary}[theorem]{Corollary}
\newtheorem{assumption}[theorem]{Assumption}

\newcommand{\norm}[1]{\left\|{#1}\right\|}

\def \caratt {{\mathds{1}}}
\def\bn{{\bar{n}}}

\def\eps{\varepsilon}

\def\e{\varepsilon}

\def\a{\alpha}

\def\p{\partial}
\def\d{\delta}
\def\g{\gamma}
\def\l{\lambda}

\def\[{{\Big[}}
\def\]{{\Big]}}
\def\<{{\langle}}
\def\>{{\rangle}}
\def\({{\big(}}
\def\){{\big)}}

\def\bx{{\mathbf{x}}}

\def\min{{\mathord{{\rm min}}}}

\def\bb2{{\boldsymbol{2}}}

\def\={&\!\!=\!\!&}

\def\cC{{\mathcal C}}

\def\hD{\hat{D}}

\def\1{{\mathbf{1}}}

\def\geq{\geqslant}
\def\leq{\leqslant}
\def\ge{\geqslant}
\def\le{\leqslant}

\def\eps{\varepsilon}

\def\a{\alpha}

\def\p{\partial}
\def\g{\gamma}
\def\l{\lambda}

\def\[{{\Big[}}
\def\]{{\Big]}}
\def\<{{\langle}}
\def\>{{\rangle}}

\def\bx{{\mathbf{x}}}

\def\min{{\mathord{{\rm min}}}}

\def\={&\!\!=\!\!&}
\def\bt{\begin{theorem}}
\def\et{\end{theorem}}
\def\bl{\begin{lemma}}
\def\el{\end{lemma}}
\def\br{\begin{remark}}
\def\er{\end{remark}}
\def\bx{\begin{Examples}}
\def\ex{\end{Examples}}
\def\bd{\begin{definition}}
\def\ed{\end{definition}}
\def\bp{\begin{proposition}}
\def\ep{\end{proposition}}
\def\bc{\begin{corollary}}
\def\ec{\end{corollary}}

\DeclareMathOperator*{\esssup}{ess\,sup}

\def\geq{\geqslant}
\def\leq{\leqslant}
\def\ge{\geqslant}
\def\le{\leqslant}

 \def\R{\mathbb R}
 \def\R{\mathbb R}    
\def\N{\mathbb N} 

\def\m{\mu}
  
\def\bbB{\langle\b\rangle_B}

\def\rr{r}
\def\<{\langle} \def\>{\rangle}

\def\phi{{\varphi}}

\def\0{{\mathbf{0}}}

\def\G{{\Gamma}}

\def\K{{\mathcal{K}}}

\def\b{\beta}

\allowdisplaybreaks

\begin{document}
\title{Sobolev embeddings for kinetic Fokker-Planck equations}

\author{Andrea Pascucci\thanks{Dipartimento di Matematica,
Universit\`a di Bologna, Bologna, Italy. {\bf Orcid:} 0000-0001-8837-5568 \textbf{e-mail}:
andrea.pascucci@unibo.it} \and Antonello Pesce\thanks{Dipartimento di Matematica, Universit\`a di
Bologna, Bologna, Italy. {\bf Orcid:} 0000-0002-1641-852X \textbf{e-mail}:
antonello.pesce2@unibo.it}}

\date{}

\maketitle

\vspace{-15pt}
\begin{center} Dedicated to Ermanno Lanconelli on the occasion of his 80th birthday
\end{center}
\begin{abstract}
We introduce intrinsic Sobolev-Slobodeckij spaces for a class of ultra-parabolic Kolmogorov type
operators satisfying the weak H\"ormander condition. We prove continuous embeddings into Lorentz
and intrinsic H\"older spaces. We also prove approximation and interpolation inequalities by means
of an intrinsic Taylor expansion, extending analogous results for H\"older spaces. The embedding
at first order is proved by adapting a method by Luc Tartar which only exploits scaling properties
of the intrinsic quasi-norm, while for higher orders we use uniform kernel estimates.
\end{abstract}


\bigskip\noindent
{\bf Statements and Declarations:} The authors have no relevant financial or non-financial
interests to disclose.

\bigskip\noindent
{\bf Data Availability Statement:} Data sharing not applicable to this article as no datasets were
generated or analysed during the current study.

\bigskip\noindent
{{\bf Acknowledgements.} The authors are members of INDAM-GNAMPA. }
\section{Introduction}

In this paper we develop a functional framework for the study of kinetic Fokker-Planck equations.
Specifically, we introduce intrinsic Sobolev spaces suitably related to a system of H\"ormander's
vector fields: our main results are embedding, interpolation and approximation theorems that are
the basic tools in many problems concerning partial differential equations.

Let $(t,x)$ denote {a} point in $\R\times\R^{N}$ and, for fixed $d\le N$, consider the vector fields 
\begin{equation}\label{vector}
 \p_{x_{1}},\dots,\p_{x_{d}}\quad\text{ and }\quad 
 Y:=\langle Bx,\nabla_{x}\rangle+\p_{t},
\end{equation}
where $B$ is a constant $N\times N$ matrix and $\nabla_{x}=(\p_{x_{1}},\dots,\p_{x_{N}})$. We
assume the H\"ormander's condition is satisfied:
\begin{equation}\label{horm}
  \text{\rm rank Lie}(\p_{x_{1}},\dots,\p_{x_{d}},Y)=N+1.
\end{equation}
The classical example we have in mind is the Langevin kinetic model, given by the system of
stochastic differential equations
\begin{equation}\label{SDE}
  \begin{cases}
    dV_{t}=dW_{t}, \\
    dP_{t}=V_{t}dt,
  \end{cases}
\end{equation}
where $W$ is a $d$-dimensional Brownian motion. Here the processes $V$ and $P$ represent the
velocity and position of a system of $d$ particles. The forward Kolmogorov (or Fokker-Planck)
operator of \eqref{SDE}, written in terms of the variables $x=(v,p)\in\R^{d}\times\R^{d}$, is in
the form of a sum of squares of the vector fields $\p_{v_{1}},\dots,\p_{v_{d}}$ plus a drift (or
transport term) $Y_{0}$, precisely
\begin{equation}\label{Lang}
  \frac{1}{2}\sum_{i=1}^{d}\p_{v_i}^{2}-Y_{0},\qquad Y_{0}:=\langle v,\nabla_{p}\rangle+\p_{t}.
\end{equation}
In this example, $N=2d$ and
  $$B=\begin{pmatrix}
    0_{d} & 0_{d} \\
    I_{d} & 0_{d} \
  \end{pmatrix}$$
where $I_{d}$ and $0_{d}$ denote the $d\times d$ identity and null matrices respectively. Operator
\eqref{Lang} satisfies the H\"ormander's condition, is hypoelliptic and has a Gaussian fundamental
solution that is the transition density of the Markov process $(V,P)$ in \eqref{SDE}.


The literature on generalized Sobolev spaces for H\"ormander's vector fields is vast (see, for
instance, \cite{MR436223}, \cite{MR558675}, \cite{MR657581}, \cite{MR1218884},\cite{MR1762582}).
When dealing with the regularity properties of PDEs modeled on the vector fields
\eqref{vector}-\eqref{horm}, as for example the equation in divergence form
\begin{equation}\label{PDE}
  \sum_{i,j=1}^{d}\p_{x_{i}}\left(a_{ij}\p_{x_{j}}u\right)-Yu=0,
\end{equation}
it is standard to assign a formal weight to each of the vector fields, that is one for the
directions $\p_{x_{1}},\dots,\p_{x_{d}}$ of diffusion and is two for the drift $Y$; also,
consistently with the structure of the equation, $Y$ should be interpreted as a {\it second order
derivative} in intrinsic sense. As earlier noted in \cite{MR436223} among others, this fact raises
a question about the role of $Y$ in the definition of {\it first order} intrinsic Sobolev space
$W^{1,p}$: indeed, in the degenerate case $d<N$, the regularity properties of \eqref{PDE} strongly
rely on H\"ormander's condition and involve the second order derivative $Y$ in a crucial way. Many
remarkable results have been proven for weak solutions of \eqref{PDE}, defined as functions $u$
such that $Yu$ belongs {to} $L^{2}$, in addition to the minimal assumptions
$u,\p_{x_{1}}u,\dots,\p_{x_{d}}u\in L^{2}$ needed to write the equation \eqref{PDE} in the sense
of distributions: we refer for instance to \cite{MR1949176}, \cite{MR2729292}, \cite{Pascucci2},
\cite{MR4290567}, \cite{MR3778645}, \cite{MR3951695} and \cite{MR4444079}. In
\cite{PascucciPolidoro} a first $L^{2}-L^{\infty}$ estimate has been proven by using Moser's
approach; moreover, in \cite{MR3923847} a Harnack inequality for kinetic Fokker-Planck equations
with rough coefficients has been proven extending the De Giorgi-Nash-Moser theory.


In Section \ref{prelim} we introduce intrinsic Sobolev-Slobodeckij spaces for \eqref{vector},
denoted by $W^{k,p}_{B}$, where at first order (i.e. $k=1$) the vector field $Y$ appears as a
fractional derivative of order $1/2$:
{this approach is coherent with the scaling properties of the H\"ormander vector fields and
therefore seems suitable for the study of \eqref{PDE}.} In particular, we can give a natural
definition of weak solution $u$ of \eqref{PDE} in the Sobolev space $W^{1,2}_B$ without requiring
$Yu\in L^{2}$ as it is usually done in the literature: as far as we know, this is the first result
in this direction.

We mention that the use of fractional derivatives makes it difficult to prove embedding results by
means of representation formulas in terms of a parametrix, at least for $k=1$, as in
\cite{PascucciPolidoro} or \cite{MR4290567}. Indeed, for the proof of our main embedding result,
Theorem \ref{Embedd_1}, we use a remarkable method developed by Tartar \cite{MR1662313}, that is
only based on scaling arguments
and a characterization of Lorentz spaces given in Lemma \ref{lem_Lorentz}. 

In the following statement $\Q$ denotes the homogeneous dimension of $\R^{N+1}$ induced by the
vector fields \eqref{vector}, whose precise definition is given in \eqref{homdim}: to fix ideas,
$\Q=4d+2$ for the Fokker-Planck operator \eqref{Lang}.
\begin{theorem}[\bf $W^{1,p}_B$ embeddings]\label{Embedd_1}\hfill
\begin{itemize}
\item[i)] For $1 \le p<\Q$ we have
\begin{equation}\label{Embedd_1e2}
  W^{1,p}_B\subseteq L^{q,p},\qquad p\le q\le p^{\ast},
  \qquad \frac{1}{p^{\ast}}=\frac{1}{p}-\frac{1}{\Q},
\end{equation}
where $L^{p,q}$ denotes the Lorentz space. In particular, {$W^{1,p}_B\subseteq L^{q}$ for
$p\le q\le p^{\ast}$};

\item[ii)] {for $\Q<p<\infty$ we have}
\begin{equation}\label{Embedd_1e1}
 W^{1,p}_B\subseteq C_B^{0,1-\frac{\Q}{p}};
\end{equation}

\item[iii)] for $p=\Q$ we have
\begin{equation}\label{ae31}
  W^{1,\Q}_B\subseteq L^{q,\Q}\subseteq L^{q},\qquad q\ge \Q.
\end{equation}
Moreover, if $u\in W^{1,\Q}_B$ then for every $\l,\d>0$ we have
\begin{equation}\label{ae30}
  \int_{|u|>\d}e^{\l|u(z)|^{\frac{\Q}{\Q-1}}}dz<\infty.
\end{equation}
\end{itemize}
\end{theorem}
The Morrey embedding \eqref{Embedd_1e1} is given in terms of the optimal generalized H\"older
spaces $C^{k,\a}_{B}$ only recently introduced in \cite{MR3429628} together with an intrinsic
Taylor formula. Embeddings for higher order spaces $W^{k,p}_{B}$ are provided in Theorem
\ref{Embedd_k}. Remarkably, estimate \eqref{ae30} extends Trudinger's result \cite{MR0216286}.
Embedding results for Kolmogorov equations were also proved in \cite{MR2434049} and more recently
in \cite{MR4444114}.

We acknowledge that Tartar himself applied his approach to the Langevin operator \eqref{Lang}:
according to \cite{MR1662313}, Appendices II and III, he proved that for a function
$f=f(t,x,v)$ on $\R\times\R^{d}\times\R^{d}$, with $f, \nabla_{v}f, Y_{0}f:= (\p_{t}+v\cdot\nabla_{x})f\in L^{p}$, 
one can first prove the ``crude'' embedding estimate
\begin{equation}
  \|f\|_{q}\lesssim \|f\|_{p}+\|\nabla_{v}f\|_{p}+\|Y_{0}f\|_{p}  ,
\end{equation}
for some $q>p$ and then get the embeddings:
\begin{itemize}
  \item[-] in $L^{p^{**},p}$ if $1\le p<3d+1$, with $\frac{1}{p^{**}}=\frac{1}{p}-\frac{1}{3d+1}$;
  \item[-] in $L^{\infty}$ if $p>3d+1$;
  \item[-] in $L^{r}$ if $p=3d+1$, for any $p\le r<\infty$.
\end{itemize}

As a main motivation, our study is a first step in the development of a theory of generalized
Besov spaces for possible applications to {\it stochastic} partial differential equations: we
mention that recent results for stochastic kinetic equations were established in \cite{MR4355925}
and \cite{ZhangBesov}. Secondly, even for deterministic kinetic equations, our results improve the
known regularity estimates available in the literature by providing the natural functional
framework for weak solutions of kinetic Fokker-Planck equations.

\medskip The paper is structured as follows. In Section \ref{prelim} we state the precise
assumptions, introduce the intrinsic Sobolev and H\"older spaces and collect some preliminary
result concerning the geometric structure induced on $\R^{N+1}$ by the vector fields
\eqref{vector}. In Section \ref{altern1} we prove a first interpolation result, Proposition
\ref{p1}, that provides a simplified and equivalent definition of intrinsic Sobolev quasi-norm. In
Section \ref{Taylora} we show an intrinsic Taylor expansion, Theorem \ref{ap1}, for functions in
$W^{k,p}_{B}$ which extends the analogous results for intrinsic H\"older spaces proved in
\cite{MR3429628}. Crucial approximation and interpolation results, Theorems \ref{Approximation}
and \ref{t3}, are proven in Section \ref{Interpolation}. Section \ref{proofmain} contains the
proof of our main result, Theorem \ref{Embedd_1}, on the embeddings of $W^{1,p}_B$. Eventually, in
Section \ref{Embedding_section} we prove Theorem \ref{Embedd_k} on the higher order embeddings.
For reader's convenience, in the Appendix we recall some basic result about interpolation and
Lorentz spaces.

\medskip{In the context of our proofs we will often use the notation $A\lesssim B$, meaning that
$A\le c B$ for some positive constant $c$ which may depend on the quantities specified in the
corresponding statement.}

\section{Preliminaries}\label{prelim}
\subsection{Assumptions}\label{Assumptions}
We recall that H\"ormander's condition is equivalent to the well-known Kalman rank condition for
controllability of linear systems (cf., for instance, Section 9.5 in \cite{MR2791231}); also, it
was shown in \cite{lanpol} that, up to a change of basis, condition \ref{horm} is equivalent to
the following
\begin{assumption}[\bf H\"ormander's condition]\label{horm2}
The matrix $B$ takes the block-form
\begin{equation}\label{B}
  B=\begin{pmatrix}
 \ast & \ast & \cdots & \ast & \ast \\
 B_1 & \ast &\cdots& \ast & \ast \\
 0 & B_2 &\cdots& \ast& \ast \\
 \vdots & \vdots &\ddots& \vdots&\vdots \\
 0 & 0 &\cdots& B_{\rr}& \ast
  \end{pmatrix}
\end{equation}
where $B_j$ is a $(d_{j-1}\times d_j)$-matrix of rank $d_j$ with
\begin{equation}
  d\equiv d_{0}\geq d_1\geq\dots\geq d_{\rr}\geq1,\qquad \sum_{j=0}^{\rr} d_j=N.
\end{equation}
\end{assumption}
%
In general, the $\ast$-blocks in \eqref{B} are arbitrary. Our second standing assumption is the
following
\begin{assumption}[\bf Homogeneity]\label{dilat}
All the $\ast$-blocks in \eqref{B} are null.
\end{assumption}
As proven in \cite{lanpol}, Assumption \ref{dilat} is equivalent to the fact that the kinetic
Fokker-Planck operator
\begin{equation}\label{operator}
  \K:=\frac{1}{2}\sum_{i=1}^{d}\p_{x_i}^{2}-Y
\end{equation}
is homogeneous of degree two with respect to the family of dilations defined as follows: first of
all, consistently with the block decomposition \eqref{B} of $B$, we write $x\in\R^{N}$ as the
direct sum $x=x^{[0]}+\cdots+x^{[r]}$ where $x^{[i]}\in\R^{N}$ is defined as
  $$x_{k}^{[i]}=
  \begin{cases}
    x_{k} & \text{ if }\, \bar{d}_{i-1}<k\le\bar{d}_{i}, \\
    0 & \text{otherwise},
  \end{cases}\qquad \bar{d}_{i}:=\sum_{j=0}^{i}d_{j},\quad \bar{d}_{-1}:=0,\quad i=0,\dots,r.
  $$
Then, we have $\K(u(D_{\l}))=\l^{2}(\K u)(D_{\l})$ where 
\begin{equation}\label{dilations}
  D_{\l}(t,x):=(\l^{2}t,\hD_{\l}x),\qquad \hD_{\l}x:=\sum_{i=0}^{r}\l^{2i+1}x^{[i]}.
\end{equation}
For instance, the Langevin operator \eqref{Lang} is homogeneous with respect to the dilation group
$D_{\l}(t,v,p)=(\l^{2}t,\l v,\l^{3}p)$ in $\R^{2d+1}$.

\subsection{Intrinsic H\"older and Sobolev spaces}\label{SP}
In this section we recall the definition of intrinsic H\"older space as given in \cite{MR3429628}
and introduce a notion of intrinsic Sobolev space, naturally associated to the system of
vector fields \eqref{vector}.

Let $h\mapsto e^{h X }z$ denote the integral curve of a Lipschitz vector field $X$ starting from
$z\in\Rdd$, defined as the unique solution of
\begin{equation}
\begin{cases}
 \frac{d}{dh}e^{h X }z= X\left(e^{h X }z\right),\qquad  &h\in\R, \\
 e^{h X }z\vert_{h=0}= z.
\end{cases}
\end{equation}
For the vector fields in \eqref{vector}, we have
\begin{equation}\label{eq:def_curva_integrale_campo}
 e^{h \p_{x_{i}} }(t,x)=(t,x+h \mathbf{e}_i),\qquad 
 e^{h Y }(t,x)=(t+h,e^{h B}x),
\end{equation}
where $\mathbf{e}_i$ is the $i$-th element of the canonical basis of $\R^{N}$.
\begin{definition}
Let $m_{X}$ be a formal weight associated to the vector field $X$. For $\a\in\,]0,m_{X}]$, we say
that $u\in C_{X}^{\alpha}$ if the quasi-norm
\begin{equation}
 \|u\|_{C^{\a}_{X}}:=
 \sup_{z\in\Rdd\atop h\in\R\setminus\{0\}} \frac{
 \left|u\left(e^{h X }z\right)-
 u(z)\right|}{|h|^{\frac{\alpha}{m_{X}}}}
\end{equation}
is finite.
\end{definition}

Hereafter, we set {\it the formal weight of the vector fields $\p_{x_{1}},\dots,\p_{x_{d}}$ equal
to one and the formal weight of $Y$ equal to two}, which is coherent with the homogeneity of the
Fokker-Planck operator $\K$ with respect to the dilations $D_{\l}$ in \eqref{dilations}. From
\cite{MR3429628} we recall the following
\begin{definition}[\bf Intrinsic H\"older spaces]\label{def:C_alpha_spaces}
For $\a\in\,]0,1]$ we define the H\"older quasi-norms
\begin{align}
  \norm{u}_{C^{0,\a}_{B}}&:=\sup_{\R^{N+1}}|u|+\sum_{i=1}^{d}\norm{u}_{C^{\a}_{\partial_{x_i}}}+\norm{u}_{C^{\a}_{Y}},\\
  \norm{u}_{C^{1,\a}_{B}}&:=\sup_{\R^{N+1}}|u|+\norm{\nabla_{d}u}_{C^{0,\a}_{B}}+\norm{u}_{C^{\a+1}_{Y}},
\intertext{where $\nabla_{d}:=(\p_{x_{1}},\dots,\p_{x_{d}})$ and inductively, for $n\ge 2$,}
  \norm{u}_{C^{n,\a}_{B}}&:=\sup_{\R^{N+1}}|u|+\norm{\nabla_{d}u}_{C^{n-1,\a}_{B}}+\norm{Y u}_{C^{n-2,\a}_{B}}.
\end{align}
\end{definition}

Next we introduce the intrinsic Sobolev spaces. First, as in \cite{MR1762582}, {for any $u\in
L^{p}$, with $p\ge1$,} we define the fractional {Sobolev-Slobodeckij} quasi-norm of order $s\in
\,]0,1[$ along a Lipschitz vector field $X$ as
\begin{align}
  [u]_{X,s,p}:=
  \left(\int_{\R^{N+1}}dz\int_{{|h|\le 1}}\frac{|u(e^{hX}z)-u(z)|^p}{|h|^{ps+1}}dh\right)^{\frac
  1p}.
\end{align}
\begin{definition}
For $p\ge 1$ we set
\begin{align}
  |u|_{1,p,B}&:=\|\nabla_{d}u\|_{p}+[u]_{Y,\frac 12,p},\\ \label{W2pp2}
  |u|_{2,p,B}&:=|\nabla_{d}u|_{1,p,B}+\|Yu\|_{p},
\intertext{and inductively, for 
$n\ge 3$,}\label{Wnppn}
  |u|_{n,p,B}&:=|\nabla_{d}u|_{n-1,p,B}+|Yu|_{n-2,p,B}.
\end{align}
\end{definition}

\begin{definition}[\bf Intrinsic Sobolev spaces]\label{def1}
For $p\ge 1$ we define the Sobolev quasi-norms
\begin{align}
  \|u\|_{W^{1,p}_B}&:=\|u\|_p+|u|_{1,p,B},\\
  \| u\|_{W^{2,p}_B}&:=\|u\|_p+{\|\nabla_{d}u\|_{W^{1,p}_B}}+\|Yu\|_{p},
\intertext{and inductively, for $n\ge 3$,}
  \|u\|_{W^{n,p}_B}&:=\|u\|_p+{\|\nabla_{d}u\|_{W^{n-1,p}_B}}+\|Yu\|_{W^{n-2,p}_B}.
\end{align}
\end{definition}
The following alternative definition of Sobolev quasi-norm is sometimes useful.
\begin{definition}\label{def2}
For $n\in\N$ and $p\ge 1$ we set
\begin{align}\label{W2pp}
 \vvvert u\vvvert_{W^{n,p}_B}&:=\|u\|_p+|u|_{n,p,B}. 
\end{align}
\end{definition}
Clearly we have $\|\cdot\|_{W^{n,p}_B}\ge\vvvert\cdot \vvvert_{W^{n,p}_B}$. In Section
\ref{altern1}, Proposition \ref{p1}, we prove that $\|\cdot\|_{W^{n,p}_B}$ and $\vvvert\cdot
\vvvert_{W^{n,p}_B}$ are equivalent and therefore define the same functional spaces. This means
that the intermediate orders quasi-norms are not needed to characterize $W^{n,p}_B$.

\begin{remark}
Let
 $$\lfloor u\rfloor_{Y,s,p}:=\left(\int_{\R^{N+1}}dz\int_{\R}\frac{|u(e^{h{Y}}z)-u(z)|^p}{|h|^{ps+1}}dh\right)^{\frac{1}{p}}$$
and notice that
  $$[u]_{Y,s,p}\le \lfloor u\rfloor_{Y,s,p}\le [u]_{Y,s,p}+c_{p,s}\|u\|_{p},\qquad
  c_{p,s}:=\left(\int_{|h|>1}\frac{2}{|h|^{1+ps}}dh\right)^{\frac{1}{p}}.$$
Then, if we replace $[u]_{Y,s,p}$ by $\lfloor u\rfloor_{Y,s,p}$ in Definition \ref{def1}, we get
equivalent norms.
\end{remark}

\subsection{Dilation and translation groups}\label{diltra}
Besides the homogeneity with respect to $D_{\l}$ in \eqref{dilations}, operator $\K$ in
\eqref{operator} has also the remarkable property of being invariant with respect to 
the left translations in the group law 
\begin{align}
  (t,x)\circ(s,\x)=(t+s,e^{sB}x+\x),\qquad (t,x),(s,\x)\in\R^{N+1}.
\end{align}
Indeed, a simple computation shows that, for any $z,\z\in \R^{N+1}$,
\begin{equation}\label{invariance_law}
 \z^{-1}\circ e^{\delta Y}z=e^{\delta Y}(\z^{-1}\circ z), \qquad \z^{-1}\circ e^{\delta \p_{x_i}}z=e^{\delta \p_{x_i}}(\z^{-1}\circ z),
 \qquad i=1,\dots , d,
\end{equation}
where $(t,x)^{-1}=(-t,-e^{-tB}x)$. 
Analogously, we have  (see, for instance, \cite{P22})
\begin{equation}\label{homogeneity}
 D_\l e^{\delta Y}(z)=e^{\delta\lambda^2Y}\left(D_\l z\right), \qquad D_\l e^{\delta
 \p_{x_i}}(z)=e^{\delta \l^{2j+1} \p_{x_i}}D_\l z, \qquad i=\bar{d}_{j-1}+1,\dots ,\bar{d}_{j}.
\end{equation}
{A $D_{\l}$-homogeneous norm on $\R^{N+1}$ is defined as
\begin{equation}\label{Bnorm}
  \|(t,x)\|_{B}=|t|^{\frac{1}{2}}+
  |x|_{B},\qquad |x|_{B}=\sum_{i=0}^{r}|x^{[i]}|^{\frac{1}{2i+1}},
\end{equation}}
and
\begin{equation}\label{homdim}
  \Q:= 2+\sum_{k=0}^{r}(2k+1)d_{k}
\end{equation}
is usually called the homogeneous dimension of $\R^{N+1}$ with respect to $D_{\l}$.
\begin{lemma}[{\cite{MR1751429},  Proposition $5.1$}] There exists $m=m(B)\ge 1$ such that
\begin{equation}\label{normcontrol}
 \|\z^{-1}\circ z\|_B\le m(\|\z\|_B+\|z\|_B)\quad m^{-1}\|z\|_B\le\|z^{-1}\|_B\le  m\|z\|_B \quad
 z,\z\in \R^{N+1}.
\end{equation}
\end{lemma}

\begin{remark}
Since $e^{\delta Y}z=z\circ (\delta,0)$, by \eqref{normcontrol} we have
\begin{equation}\label{normcontrol2}
  \frac{1-mc}{{m}}\|z\|_{B}\le \|e^{\delta Y}z\|_B\le m(1+c)\|z\|_B,
\end{equation}
for any $|\d|^{\frac 12}\le c\|z\|_B$ with $c\in\,]0,\frac{1}{m}[$.
\end{remark}

{\begin{remark}\label{r1} The matrix $B$ is nilpotent of degree $r+1$. In particular, {for any
$n\le r$ we have}
\begin{equation}\label{Bn}
B^n=\begin{pmatrix} 0_{\bar{d}_{n-1}\times d_0}& 0_{\bar{d}_{n-1}\times d_1} &\cdots &
0_{\bar{d}_{n-1}\times d_{r-n}}& 0_{\bar{d}_{n-1}\times (\bar{d}_{r}-\bar{d}_{r-n})}\\
\prod\limits_{j=1}^n{B}_j &0_{d_n\times d_1} &\cdots & 0_{d_n\times d_{r-n}}& 0_{d_n\times
(\bar{d}_{r}-\bar{d}_{r-n})}\\ 0_{d_{n+1}\times d_0}& \prod\limits_{j=2}^{n+1}{B}_j &\cdots &
0_{d_{n+1}\times d_{r-n}}& 0_{d_{n+1}\times (\bar{d}_{r}-\bar{d}_{r-n})}\\ \vdots &\vdots & \ddots
&\vdots &\vdots \\ 0_{d_r\times d_0}& 0_{d_r\times d_1} &\cdots &
\prod\limits_{j={r-n+1}}^{{r}}{B}_j & 0_{d_r\times (\bar{d}_{r}-\bar{d}_{r-n})}
\end{pmatrix}, 
\end{equation}
{where
  $$\prod_{j=1}^{n}B_{j}=B_{n}B_{n-1}\cdots B_{1},$$
}and ${B}^n=0$ for $n>r$. Thus
 $$e^{\delta B}=I_N+\sum_{j=1}^r\frac{B^j}{j!}\delta^j$$
is a lower triangular matrix with diagonal $(1,\dots, 1)$ and therefore it has determinant equal
to 1.
\end{remark}
}
\begin{lemma}For any $n\in \N$ {and $u\in L^{p}$, with $p\ge 1$,} we have
\begin{equation}\label{scaling}
 \|u( D_\l)\|_p=\l^{-\frac{\Q}{p}}\|u\|_p,\qquad |u( D_\l)|_{n,p,B}=\l^{n-\frac{\Q}{p}}|u|_{n,p,B}
\end{equation}
\end{lemma}
\proof The first equality follows by a simple change of variable. Next, for $i=1,\dots,d$ we have,
by \eqref{homogeneity}
\begin{align}
 \|\p_{x_i}u(D_\l)\|_p^p&=\int_{\R^{N+1}}|\p_{x_i}u(D_\l z)|^pdz=\int_{\R^{N+1}}|\l
 (\p_{x_i}u)(D_\l z)|^pdz=
\intertext{(by the change of variable $z'=D_\l z$)}
 &=\l^p\int_{\R^{N+1}}|(\p_{x_i}u)(z')|^p\l^{-\Q}dz'=\l^{p-\Q}\|\p_{x_i}u\|^p_p.\label{sc1}
\end{align}
Similarly
\begin{align}
 [u(D_\l)]_{Y,\frac 12,p}^p&=\int_{\R^{N+1}}\int_{\R}{|u\left(D_\l(e^{hY}z)\right)-u(D_\l
 z)|^p}\frac{dh}{|h|^{1+\frac{p}{2}}}dz\\
 &=\int_{\R^{N+1}}\int_{\R}{\l^{2+p}|u(e^{h\l^2Y}D_\l z)-u(D_\l z)|^p}\frac{dh}{|\l^2
 h|^{1+\frac{p}{2}}}dz=
\intertext{(by the change of variables $(h',z')=(\l^2 h, D_\l z)$)}
 &=\l^{p}\int_{\R^{N+1}}\int_{\R}{|u(e^{h'Y}z')-u(z')|^p}\frac{dh'}{|h'|^{1+\frac{p}{2}}}\l^{-\Q}dz=\l^{p-\Q}[u]_{Y,\frac
 12,p}^p.\label{sc2}
\end{align}
\eqref{sc1} and \eqref{sc2} give the second equality for $n=1$. The case $n=2$ is analogous and
{the general case $n>2$ follows by induction}.
\endproof





\section{Alternative Sobolev norms and a first interpolation result}\label{altern1}
\subsection{Intrinsic weak derivatives in $W^{n,p}_B$}
By definition, the quasi-norm $|\cdot|_{n,p,B}$ only controls weak derivatives of order $n$ and
$n-1$, which are made up of compositions of $\p_{x_1},\dots,\p_{x_d}$ and $Y$ for any possible
permutation. We show that actually a function $u\in W_{B}^{n,p}$ supports all the weak derivatives
of intrinsic order $l$, $l\le n$, and for $k\in \N_0$ and $\b\in \N_0^{N}$ we have
\begin{equation}\label{intrinsic_D}
 Y^k\p^\b u\in W^{n-l,p}_B,\quad 2k+\bbB=l,
\end{equation}
where
\begin{equation}\label{B_lenght}
 \p^\b=\p^{\b_1}_{x_1}\cdots \p^{\b_N}_{x_N}, \qquad 
 {\bbB:=\sum_{i=0}^r(2i+1)\sum_{k=1+\bar{d}_{i-1}}^{\bar{d}_{i}}\b_{k}.}
\end{equation}
Indeed these derivatives can be recovered by taking appropriate iterated commutators of the vector
fields {$\p_{x_1},\dots,\p_{x_d}$ and $Y$}: exploiting these commutators, we can also rearrange
the terms appearing in $|\cdot|_{n,p,B}$ and provide a more explicit characterization which only
make use of the intrinsic derivatives in the form \eqref{intrinsic_D}.

\medskip
\noindent First we recall some preliminary notions from \cite{MR3429628}, Section 4. 
By the structure of the matrix $B$, for any $n=0,\dots,r$ and $v\in \R^N$ we have
\begin{equation}\label{eqa2bbb}
  B^nv \in \bigoplus_{i=n}^rV_i, \qquad V_i:=\{x^{[i]}\mid x\in\R^N\},
\end{equation}
and $B^n=0$ for $n>r$. In particular, if $v\in V_0$ then we have
\begin{equation}
  B^nv\in V_n, \qquad n=0,\dots, r.
\end{equation}
Moreover there exist subspaces $$V_{0,r}\subseteq V_{0,r-1}\subseteq \cdots \subseteq
V_{0,1}\subseteq V_{0,0}:=V_0 $$ such that the linear map
\begin{equation}\label{eqa2bb}
  \psi_n: V_{0,n}\longrightarrow V_n,\qquad \psi_n(v):=B^nv
\end{equation}
is bijective. For $v\in V_0$, we introduce the following iterated commutators
  $${X^{(0)}_v:=
  \sum_{k=1}^dv_k\p_{x_k},}$$
and recursively
\begin{equation}\label{Commutator}
  {X^{(n)}_v:=[X^{(n-1)}_v,Y]=X^{(n-1)}_vY-YX^{(n-1)}_v, \qquad n\in \N.}
\end{equation}

\begin{lemma}
Let $u\in W^{n,p}_B$. Then, for any $i\in \N_{0}$, $2i+1\le n$, and $v\in V_0$ we have
\begin{equation}\label{eqa1}
  X^{(i)}_vu\in W^{n-2i-1,p}_B.
\end{equation}
\end{lemma}
\proof We use an induction argument on $n$. If $n\le 2$ there is nothing to prove because
$X^{(0)}_v$ is a linear combination of the vector fields $\p_{x_1},\dots,\p_{x_{d}}$ and the
thesis follows by definition.

Assume \eqref{eqa1} is true for some fixed $n>2$ and let us prove it for $n+1$. We proceed by
induction on $i$. For $i=0$, again there is nothing to prove. We assume \eqref{eqa1} for some
$i>0$ such that $2(i+1)+1\le n +1$ and prove it for $i +1$. If $u\in W^{n+1,p}_B$ then
  $$X^{(i+1)}_vu=X^{(i)}_vYu-YX^{(i)}_vu.$$
Here $Yu\in W^{n-1,p}_B$ by definition and therefore $X^{(i)}_vYu\in W^{n-1-(2i+1),p}_B$ by the
inductive hypothesis on $n$. On the other hand $X^{(i)}_vu\in W^{n+1-(2i+1),p}_B$ by the inductive
hypothesis on $i$ and therefore $YX^{(i)}_vu\in W^{n+1-(2i+1)-2,p}_B$ by definition. Then
\eqref{eqa1} holds for $n+1$, for any $i$, $2i+1\le n+1$, and this concludes the proof.
\endproof

\begin{proposition}\label{P01}Let $u\in W^{n,p}_B$. Then for any $i\in \N_{0}$, $2i+1\le n$, we
have
\begin{equation}\label{eqa4}
  {\p_{x_j}u\in W^{n-2i-1,p}_B, \qquad j=1+\bar{d}_{i-1},\dots, \bar{d}_{i}.}
\end{equation}
\end{proposition}
\proof By induction it is not difficult to prove that
\begin{equation}\label{eqa2}
 X^{(i)}_v\phi=\langle B^iv,\nabla_{x} \phi\rangle,\qquad \phi \in C^{\infty}.
\end{equation}
Next, since $\psi_i$ in \eqref{eqa2bb} is bijective, for every $j=1+\bar{d}_{i-1},\dots,
\bar{d}_{i}$ there exists $w_j\in V_{0,i}$ such that $B^iw_j=\mathbf{e}_j\in V_i$. Then
$f:=X^{(i)}_{w_j}u\in W^{n-2i-1,p}_B$ is such that $$\int_{\R^{N+1}}f(z)\phi(z)dz=-
\int_{\R^{N+1}}u(z)X^{(i)}_{w_j}\phi(z)dz=- \int_{\R^{N+1}}u(z)\p_{x_j}\phi(z)dz,\quad \phi\in
C_0^{\infty},$$ which means that $f$ {is the weak derivative $\p_{x_j}u$}.
\endproof

By Proposition \ref{P01} and the definition of intrinsic Sobolev spaces we eventually infer the
following:
\begin{corollary}\label{P02}Let $u\in W^{n,p}_B$. Then, for any $k\in \N_0$, $\b\in \N_0^{\b}$ such that $2k+\bbB=l\le
n$, we have
\begin{equation}\label{eqa3}
  Y^k\p^\b u\in W^{n-l,p}_B.
\end{equation}
\end{corollary}

{\begin{corollary} The following quasi-norms are equivalent:
\begin{itemize}
\item[i)] $|u|_{n,p,B}$;
\item[ii)]
\begin{equation}\label{eqa5}
 \sum_{2k+\bbB=n}\|Y^k\p^\b u\|_p+\sum_{2k+\bbB=n-1}[Y^k\p^\b u]_{Y,\frac 12,p};
\end{equation}
\item[iii)]
\begin{equation}\label{eqa6}
 \begin{cases}
  \sum\limits_{2k+\bbB=n-1}|Y^k\p^\b u|_{1,p,B}, & n=2l+1, \ l\in \N, \\
  \sum\limits_{2k+\bbB=n-1}|Y^k\p^\b u|_{1,p,B}+\|Y^{l}u\|_p, &n=2l, \ l\in \N.
 \end{cases}
\end{equation}
\end{itemize}
\end{corollary}}
\proof By induction it is not difficult to check that $|u|_{n,p,B}$ controls all the $L^{p}$-norms
of the $n$th-order derivatives that are compositions of $Y$ and $\p_{x_1},\dots,\p_{x_d}$ for any
possible permutation, as well as the fractional quasi-norms of the $(n-1)$th-order derivatives.
Then it suffices to note that, by \eqref{Commutator} we have
\begin{equation}\label{Commutator2}
  X^{(0)}_vY^n=\sum_{i=0}^n \binom{n}{i}Y^iX_v^{(n-i)}, \qquad v\in V_0;
\end{equation}
then, proceeding as in the proof of Proposition \ref{P01} to rearrange the derivatives, we get
\eqref{eqa5}.

Moreover, by definition we have
 $$\sum_{2k+\bbB=n-1}|Y^k\p^\b u|_{1,p,B}=\sum_{2k+\bbB=n-1}\left(\sum_{j=1}^d\|{\p_{x_{j}}Y^k\p^\b} u\|_{p}+[Y^k\p^\b u]_{Y,\frac12,p}\right),$$
and the fractional part of the quasi-norm coincides with \eqref{eqa5}. As for the first term in
the sum we need to distinguish two cases: if $n=2l+1$ for some $l\in \N$, then we get an
equivalence with \eqref{eqa5} by rearranging the derivatives as in the proof of Proposition
\ref{P01}; indeed, compared to $W^{n-1,p}_B$ we have the additional set of Euclidean derivatives
$\p_{x_j}$, $j=1+\bar{d}_{l-1},\dots, \bar{d}_{l}$ which can be recovered from
$\sum_{j=1}^d\p_{x_j}Y^lu$ by \eqref{Commutator2}, and similarly for the mixed derivatives. If
$n=2l$ we have the derivatives $Y^lu$ that cannot be written as sums of iterated commutators and
thus we get \eqref{eqa6}.
\endproof

\subsection{Interpolation inequality and equivalence of the norms $\|\cdot\|_{W^{n,p}_B}$ and $\vvvert\cdot \vvvert_{W^{n,p}_B}$}
\begin{proposition}\label{p1} Let $1\le n{< m}$ and $p\ge1$. {There exists $c=c(m,p,B)$ such that
\begin{equation}\label{interp_1}
  |u|_{n,p,B}\le c \left(\eps |u|_{m,p,B}+\eps^{-\frac{n}{m-n}}\|u\|_p\right), \qquad u\in W^{m,p}_{B}, \quad \e>0.
\end{equation}}
In particular the norms $\|\cdot\|_{W^{n,p}_B}$ and $\vvvert \cdot\vvvert_{W^{n,p}_B}$ are
equivalent.
\end{proposition}
\proof The proof is based on a two-step induction.

\medskip\noindent \textit{Step 1: case $n=1$ and $m=2$}. The estimate
 $$\|\p_{x_i}u\|_p\lesssim \|u\|_p+\|\p_{x_ix_i}u\|_p,\qquad i=1,\dots,d,$$
is standard (cf. for instance \cite{MR2424078}, Chapter 5). On the other hand,
by Fubini's Theorem we have
\begin{align}
 [u]^p_{{Y,\frac 12,p}}=\int_{-1}^1\frac{dh}{|h|^{1+\frac p2}}\int_{\R^{N+1}}J_p(z,h)dz,\qquad
 J_p(z,h):=|u(e^{hY}z)-u(z)|^p.
\end{align}
By the mean value theorem along the vector field $Y$, for every $z\in\R^{N+1}$ and $h\in [-1,1]$,
$h\neq 0$, there exists $|\bar h|\le |h|$ such that $|J_p(z,h)|\le |Yu(e^{\bar{h} Y}z)|^p|h|^p$:
then, by a change of variable  and recalling Remark \ref{r1},
  $$ [u]^p_{{Y,\frac 12,p}} \le 2\|Yu\|^p_p \int_{0}^1 |h|^{\frac p2-1} dh \le \frac{4}{p}\|Yu\|^p_p.$$
Thus we obtain
\begin{equation}\label{interp_2}
  |u|_{1,p,B}\lesssim |u|_{2,p,B}+\|u\|_p.
\end{equation}
The thesis follows by a scaling argument: indeed, applying \eqref{interp_2} to
$u(D_{\eps^{-1}}\cdot)$, by \eqref{scaling} we get
  $$\eps^{1-\frac{\Q}{p}}|u|_{1,p,B}\lesssim \eps^{2-\frac{\Q}{p}}|u|_{2,p,B}+\eps^{-\frac{\Q}{p}}\|u\|_{p}.$$

\medskip
\noindent \textit{Step 2: induction on $n,m$ with $m=n+1$.} {We first prove the preliminary
interpolation inequality: 
\begin{equation}\label{interp_2bis}
  \|Yu\|_p\lesssim \e [Yu]_{Y,\frac 12,p}+\e^{-1} [u]_{Y,\frac 12,p}, \qquad u\in W^{3,p}_{B}, \quad \e>0.
\end{equation}
We have
  $$u(e^{Y}z)-u(z)-Yu(z)=\int_0^1\left(Yu(e^{\d Y}z)-Yu(z)\right)d\delta,$$
and therefore
\begin{align}
  \|Yu\|^p_p&=\int_{{\R^{N+1}}}\left|u(e^Yz)-u(z)-\int_0^1\left(Yu(e^{\d
  Y}z)-Yu(z)\right)d\delta\right|^pdz\le
\intertext{(by the triangular and H\"older inequalities)}
  &\lesssim \int_{\R^{N+1}}|u(e^{Y}z)-u(z)|^pdz+\int_{\R^{N+1}}\int_0^1|Yu(e^{\delta Y}z)-Yu(z)|^pd\delta dz=:I_1+I_2,
\end{align}
where
\begin{equation}
 I_2=\int_{\R^{N+1}}\int_0^1\frac{|Yu(e^{\delta Y}z)-Yu(z)|^p}{\delta^{1+\frac p2}}\delta^{1+\frac
 p2}d\delta dz\le \frac 12 [Yu]_{Y,\frac 12,p}
\end{equation}
and
\begin{align}
 I_1&\le \int_{\R^{N+1}}\int_0^1|u(e^{Y}z)-u(e^{\delta Y}z)|^pd\delta
 dz+\int_{\R^{N+1}}\int_0^1|u(e^{\delta Y}z)-u(z)|^pd\delta dz=
\intertext{(by the change of variables $z'=e^{Y}z$ and $\d'=\d-1$)}
 &=\int_{\R^{N+1}}\int_{-1}^0|u(e^{\d'Y}z')-u(z')|^pd\delta'
 dz'+\int_{\R^{N+1}}\int_0^1|u(e^{\delta Y}z)-u(z)|^pd\delta dz\le [u]_{Y,\frac 12,p},
\end{align}
reasoning as for $I_2$ in the last step. Then \eqref{interp_2bis} follows by a scaling argument.}

\medskip Next we prove that if, for some $\bn\in\N$, \eqref{interp_1}
holds with $n=\bn$, $m=\bn+1$ 
then it also holds with $n=\bar{n}+1$ and $m=\bar{n}+2$. By Step 1, \eqref{eqa6} and
\eqref{interp_2bis}, if $\bar n+1$ is \textit{even} we have
\begin{align}
 |u|_{\bn+1,p,B}&{\lesssim}\sum_{2k+\bbB=\bn}|Y^{k}\p^{\b}u|_{1,p,B}+{\|Y^{\frac{\bn+1}{2}}u\|_p}\\ &\le c_1\sum_{2k+\bbB=\bn}\left(\eps\left(
 |Y^{k}\p^{\b}u|_{2,p,B}+[Y^{\frac{\bn+1}{2}}u]_{Y,\frac 12,p}\right)+\eps^{-1}\left(\|Y^{k}\p^{\b}u\|_{p}+
 [Y^{\frac{\bn-1}{2}}u]_{Y,\frac 12,p}\right)\right) \\
 &\le c_1\left(\eps
 |u|_{\bn+2,p,B}+\eps^{-1} |u|_{\bn,p,B}\right)\le
\intertext{(by the inductive hypothesis)}
 &\le c_1\eps |u|_{\bn+2,p,B}+c_1c_2\eps^{-1}\eps_1|u|_{\bn+1,p,B}+c_1c_2\eps^{-1}\eps^{-\bn}_1\|u\|_p.
\end{align}
If $\bar n+1$ is \textit{odd}, by \eqref{eqa6} we derive the same estimate only exploiting Step 1.
To conclude it suffices to take $\eps_1=\frac{\e}{2c_1c_2}$. 

\medskip
\noindent\textit{Step 3: backward induction on $n$.} Let $m\in\N$, $m>2$, be fixed. We prove that
if \eqref{interp_1} is true for $m$ and $n=\bn$ for some $\bar{n}\in\{2,\dots,m-1\}$, then it is
also true for $m$ and $n=\bar{n}-1$.

By Step 2 we have
\begin{align}
 |u|_{\bn-1,p,B}&\lesssim \eps_{1} |u|_{\bn,p,B}+\eps_{1}^{-(\bn-1)}\|u\|_p\lesssim
\intertext{(by the inductive hypothesis)}
 &\lesssim \eps_{1}\left(\eps_{2}|u|_{m,p,B}+\eps^{-{\frac{\bn}{m-\bn}}}_2\|u\|_p\right)+\eps_{1}^{-(\bn-1)}\|u\|_p.
\end{align}
Letting now $\e=\e_{1}\e_{2}$ and $\eps_{1}=\e^{\frac{1}{m-(\bn-1)}}$, we get
  $$|u|_{\bn-1,p,B}\lesssim \eps |u|_{m,p,B}+\eps^{-\frac{\bn-1}{m-(\bn-1)}}\|u\|_p,$$
which concludes the proof.
\endproof

\section{Taylor expansion in $W^{n,p}_B$}\label{Taylora}
According to \cite{MR3429628}, the \emph{$n$-th order $B$-Taylor polynomial of $u$ around
$\z=(s,\xi)$} is formally defined as
\begin{equation}\label{eq:def_Tayolor_n}
 T_n u(\z,z):= \sum_{{0\leq 2 k + \bbB \leq n}}\frac{(t-s)^k (x-e^{(t-s)B}\xi)^{\beta}}{k!\,\beta!}
 Y^k \partial_{\xi}^{\beta}u(s,\xi),\qquad
 z=(t,x)\in\R^{N+1},
\end{equation}
with $\bbB$ as in \eqref{B_lenght}.
The main result of this section is the following.
\begin{theorem}\label{ap1}
{Let $n\in\N_{0}$ and $p\ge 1$. There exists $c=c(n,p,B)$ such that, for any $u\in W^{n+1,p}_B\cap
C^{\infty}$ we have
\begin{align}
  \label{e14b}
  \|u-T_nu(\cdot\circ \z,\cdot) \|_p&\le c \|\z\|^{n+1}_B\|u\|_{W^{n+1,p}_B}, \qquad \z\in\R^{N+1}.
\end{align}}
\end{theorem}

The proof is based on an induction procedure developed in \cite{MR3429628} to derive the
$C^{n,\a}_B$ estimate of the remainder. For completeness, here we give a fairly comprehensive
presentation of the main lines, and refer to \cite{MR3429628} for the details of the construction.
To simplify the exposition we first split the proof in different steps, corresponding to
particular cases of \eqref{e14b}.

\begin{lemma}{There exists $c=c(n,p)$ such that, for any $u\in W^{n+1,p}_B\cap C^{\infty}$ and $\d\in \R$, we have}
\begin{align}
  \label{e11b}
  \|u(e^{\delta Y}\cdot)-\sum_{k=0}^{[n/2]}\frac{\d^{k}}{k!}Y^{k}u\|_p &{\le c}
  |\delta|^{\frac{n+1}{2}} \|u\|_{W^{n+1,p}_{B}},
\end{align}
where $[n/2]$ denotes the integer part of $n/2$.
\end{lemma}
\proof We first check that, for $u\in W^{1,p}_{B}$, we have
\begin{align}
 \label{e11}
 \|u(e^{\delta Y}\cdot)-u\|_p&{\lesssim}|\delta|^{\frac 12} {[u]_{Y,\frac 12,p}}, \qquad {|\delta|\le 1}.
\end{align}
Without loss of generality, we assume $\delta{\in (0,1]}$. Adding and subtracting $u(e^{h
Y}\cdot)$ and integrating on $h\in[0,\delta]$, we have
\begin{align}
 \delta \|u(e^{\delta Y}\cdot)-u\|^p_p&\lesssim \int_{0}^\delta\|u(e^{\delta Y}\cdot)-u(e^{h Y}\cdot)\|^p_p dh
 +\int_0^\delta\|u(e^{h  Y}\cdot)-u\|^p_p dh =:I_1+I_2.
\end{align}
Then we have
 $$I_2=\int_{\R^{N+1}}\int_0^\delta\frac{|u(e^{h Y}z)-u(z)|^p}{h
 ^{1+\frac{p}{2}}}h^{1+\frac{p}{2}}dh\,{dz}
 {\le\frac 12} |\delta|^{1+\frac{p}{2}}[u]^p_{Y,\frac 12,p}.$$
By the change of variable {$\bar z=e^{\delta Y}z$, $\bar h=h-\delta$} and Remark \ref{r1}, the
term $I_1$ is analogous, and thus we get \eqref{e11}. Similarly, we also see that if $u\in
W^{2,p}_{B}$ then
\begin{align}\label{e11c}
 \|u(e^{\delta Y}\cdot)-u\|_p&{\lesssim}|\delta| \|Yu\|_{p}{\lesssim}|\delta| \|u\|_{W^{2,p}_B}.
\end{align}
More generally, for $u\in W^{n+1,p}_B\cap C^{\infty}$, 
by the mean-value theorem along the vector field $Y$, for some $\bar{\d}$ such that $|\bar{\d}|\le
|\d|$ we have
\begin{equation}\label{e11b2}
  u(e^{\delta Y}z)-\sum_{i=0}^{[n/2]}\frac{\d^{i}}{i!}Y^{i}u(z)=
  \frac{\d^{[n/2]}}{[n/2]!}\left(Y^{[n/2]}u(e^{\bar{\d}Y}z)-Y^{[n/2]}u(z)\right).
\end{equation}
Now, if $n=2h$ for some $h\in\N$ then $Y^{h}u\in W^{1,p}_{B}$: thus \eqref{e11b} follows by
combining \eqref{e11b2} with \eqref{e11} applied to $Y^{h}u$. Similarly, if $n=2h+1$ for some
$h\in\N$ then \eqref{e11b} follows by combining \eqref{e11b2} with \eqref{e11c} applied to
$Y^{h}u\in W^{2,p}_{B}$.
\endproof

\begin{lemma}\label{la1}
{There exists $c=c(p,B)$ such that
\begin{equation}\label{e10bis}
  \|u(\cdot\circ (0,\x))-u\|_p {\le c} |\x|_B\|u\|_{W^{1,p}_B}, \qquad u\in W^{1,p}_B, \quad \x\in \R^{N}.
\end{equation}}
\end{lemma}
\proof We have the standard inequality
\begin{equation}\label{e10}
 \|u(e^{\delta \partial_{x_{i}}}\cdot)-u\|_p{\le}|\delta|\|\partial_{x_{i}}u\|_{p}, \qquad
 i=1,\dots,d.
\end{equation}
Notice that, for $u\in W^{1,p}_{B}$, we have $L^{p}$-bounds only on the first $d$ spatial derivatives. 
Thus, in order to prove \eqref{e10bis} we must exploit estimate \eqref{e11} and connect any
arbitrary point $z=(t,x)\in\R^{N+1}$ to $z\circ (0,\x)=(t,x+\x)$ through a chain of integral
curves associated only to the vector fields $\p_{x_1},\dots,\p_{x_{d}}$ and $Y$.

To do so, we define a sequence of points 
{$(z_k=(t,{x(k)}))_{k=0,\cdots, r}$} adjusting, at any step $k$, the set of variables of {the
layer $V_{k}$} in \eqref{eqa2bbb}. Following \cite{MR3429628}, Lemma 4.22 we set
  $$v_0=\frac{\x^{[0]}}{|\x^{[0]}|}, \qquad |\delta_0|=|\x^{[0]}|,$$ and
  $$z_{-1}=z, \qquad z_0=\gamma^{(0)}_{v_0,\delta_0}(z_{-1}):=e^{\d_0 \langle
  v_0,{\nabla}\rangle}z=(t,{x+\x^{[0]}}).$$
For $k=1,\dots,r$ let
 $$z_k=\g^{(k)}_{v_k,\delta_k}(z_{k-1}):=e^{-\d_k^2
 Y}\left(\g^{(k-1)}_{v_k,-\delta_k}\left(e^{\d_k^2
 Y}\left(\g^{(k-1)}_{v_k,\delta_k}(z_{k-1})\right)\right)\right), \qquad
 \delta_k=|w_k|^{\frac{1}{2k+1}},$$
{where $v_k=w_k/|w_k|$ and $w_k$ is the unique vector in $V_{0,k}\subset V_0$ such that
{$B^kw_k=\x^{[k]}+x^{[k]}-x^{[k]}(k-1)$.} Importantly, it can be proven by induction that,
for any $v\in V_0$ we have 
  $$\g^{(k)}_{v,\d}(t,x)=(t,x+S_k(\d)v), \qquad S_{k}(\d)v=(-1)^k\sum_{h\in\N^k\atop |h|\leq
  r}\frac{(-B)^{|h|}}{h!}\d^{2|h|+1}v\in \bigoplus_{j=k}^r V_j.$$
In other words the flow $\g^{(k)}_{v,\d}$ only affects the set of variables $[k:r]$. Moreover,
$\g^{(r)}_{v_r,\d_r}(z_{r-1})=(t,x+\x)$ by construction. Notice also that at any step, $\d_k$ does
not depend on $x$ and also, the specific choice of $w_k$ implies $\d_k\le c_B |\x|_B$ (cf.
\cite{MR3429628}, Lemma 4.22).

We are ready to prove \eqref{e10bis}: by the Minkowski inequality we have $$\|u(\cdot\circ
(0,\x))-u\|^p_p\lesssim \sum_{k=0}^r\int_{\R^{N+1}}|u(z_k)-u(z_{k-1})|^p dz =
\sum_{k=0}^r\|u(\g_{v_k,\d_k}^{(k)})-u\|_{p}^{{p}},$$ where we exploited the changes of variables
$z'=z_{k-1}\equiv z_{k-1}(z)$ {in the last step: here we use the fact that, by
Remark \ref{r1}, the Jacobian of the change of variables has determinant equal to one.} 
Then the proof is completed once we have proved that, {for any $k, i\in\{0,\dots,r\}$ we have}
\begin{equation}\label{e13}
{\|u(\g_{v_i,\d_i}^{(k)})-u\|_p\lesssim \d_i \|u\|_{W^{1,p}_B} \lesssim |\x|_B\|u\|_{W^{1,p}_B}.}
\end{equation}
We proceed by induction on $k$. The case $k=0$ follows from \eqref{e10}. Assume now \eqref{e13}
holds for some $k\in\{0,\dots,r-1\}$: as before we have  {
\begin{align}
\|u (\g_{v_i,\d_i}^{(k+1)})-u\|_p\lesssim &\|u(\g_{v_i,\d_i}^{({k})})-u\|_p +\|u(e^{\d^2_{i}
Y}\cdot)-u\|_p \\ & +
 \|u(\g_{v_i,-\d_i}^{({k})})-u\|_p+\|u(e^{-\d^2_i Y}\cdot)-u\|_p,
\end{align}}
and the thesis follows from the inductive step and \eqref{e11}.
\endproof

We are ready to prove Theorem \ref{ap1}.

\proof[Proof of Theorem \ref{ap1}] We prove that, for any $n\in\N_{0}$, $\z\in\R^{N+1}$ and $u\in
W^{n+1,p}_B\cap C^{\infty}$, we have
\begin{equation}\label{e14}
 \|u(\cdot\circ \z)-T_{n}u(\cdot,\cdot\circ \z)\|_p {\le c} \|\z\|^{n+1}_B\|u\|_{W^{n+1,p}_B}.
\end{equation}
Estimate \eqref{e14b} follows from \eqref{e14} since
  $$\|u-T_nu(\cdot\circ \z,\cdot) \|_p=\|u(\cdot\circ \z^{-1})-T_{n}u(\cdot,\cdot\circ
  \z^{-1})\|_p.$$
Now, for $z=(t,x),\z=(s,\x)\in\R^{N+1}$, we write
\begin{align}
  u(z\circ \z)-T_{n}u(z, z\circ \z)=\LaTeXunderbrace{u(z\circ \z)-T_{n}u(e^{s Y}z,z\circ \z)}_{=:F_{1}(z,\z)}
  +\LaTeXunderbrace{T_{n}u(e^{s Y}z,z\circ \z)-T_{n}u(z,z\circ \z)}_{=:F_{2}(z,\z)}.
\end{align}
By definition
\begin{align}\label{ae12}
  {e^{s Y}z=(t+s,e^{sB}x), \qquad z\circ \z=(t+s,\x+e^{sB}x).}
\end{align}
Hence $F_{1}(z,\z)$ contains increments that only differ in the spatial variables, while
$F_{2}(z,\z)$ contains increments that only differ along $Y$.

To estimate $F_{2}(z,\z)$ we first notice that the increments in the Taylor polynomials appearing
in $F_{2}(z,\z)$ are given by
\begin{align}\label{ae17}
  (e^{s Y}z)^{-1}\circ(z\circ \z)=(0,\x), \qquad z^{-1}\circ (z\circ \z)=(s,\x),
\end{align}
thus we have
\begin{align}
 F_{2}(z,\z)&=  \sum_{\bbB\leq n}\frac{\x^{\b}}{\b!} (\partial^{\beta}u)(e^{s Y}z)
      - \sum_{2k + \bbB\leq n}\frac{s^k\x^{\b}}{k!\,\b!}(Y^k \partial^{\beta}u)(z) \\
 &=\sum_{\bbB\leq n}\frac{\x^{\b}}{\b !} \left((\partial^{\beta} u)(e^{s Y}z)
     - \sum_{2k\leq n-\bbB} \frac{s^k}{k!} (Y^k\p^{\b} u)(z)\right).
\end{align}
Taking the $L^p$ norm in $dz$ and using \eqref{e11b} {for $\partial^{\beta} u\in W^{n-\bbB}_B$ by
Corollary \ref{P02}}, we get
  $$ \|F_{2}(\cdot,\z)\|_p\lesssim \sum_{\bbB\leq n}|\x|_B^{\bbB}|s|^{\frac{n-\bbB+1}{2}}\|\p^{\b}u\|_{W^{n+1-\bbB,p}_B}\lesssim \|\z\|^{n+1}_B\|u\|_{W^{n+1,p}_B}.$$
It remains to prove
\begin{equation}\label{ae14}
  \|F_{1}(\cdot,\z)\|_p\lesssim\|\z\|^{n+1}_B\|u\|_{W^{n+1,p}_B}.
\end{equation}
First notice that, by a change of variable we have $$\|F_{1}(\cdot,\z)\|_p=\|u(\cdot\circ
(0,\x))-T_nu(\cdot,\cdot \circ (0,\x))\|_p.$$ The case $n=0$ corresponds to Lemma \ref{la1}. Next
we assume that \eqref{ae14} holds for $n=\bar n -1$ and prove it for $n=\bar n$. We have
\begin{align}
  &u(z\circ (0,\x))-T_{\bar{n}}u(z,z \circ (0,\x))\\
  &\qquad =u(t,x+\x)-T_{\bar{n}}u((t,x),(t,x+\x))\\
  &\qquad =\LaTeXunderbrace{u(t,x+\x)-T_{\bar{n}}u((t,\bar x),(t,x+\x))}_{=:F_{11}(z,\z)}+
  \LaTeXunderbrace{T_{\bar{n}}u((t,\bar
  x),(t,x+\x))-T_{\bar{n}}u((t,x),(t,x+\x))}_{=:F_{12}(z,\z)},
\end{align}
where $\bar x$ the point in $\R^{N}$ defined by
\begin{equation}\label{ae15}
  \bar{x}^{[i]}=\begin{cases} x^{[i]} & \text{if  }  {2i+1 {\le}\bar n}, \\
  x^{[i]}+\x^{[i]} &
  \text{if  }  {2i +1{>} \bar n}.
  \end{cases}
\end{equation}
Notice that $(x+\x-\bar x)^{\b}=\x^{\b}$ for $\bbB\le \bar n$ and $|x-\x-\bar x|_B\le |\x|_B$.
{For $x\in\R^{N}$, we introduce the notation
\begin{equation}\label{ae11}
  {x^{[i:j]}=\sum_{k=i}^{j}x^{[k]}\qquad 0\le i<j\le r.}
\end{equation}
} Then
\begin{equation}\label{ae16}
 T_{\bar n}u((t,\bar x),(t,x+\x))=\sum_{\bbB\le \bar n}\frac{1}{\b!}\p^{\b}u(t,{x+\x^{[[\frac{\bar n+1}{2}]:r]}})\x^{\b}.
\end{equation}
Now, since $u\in W^{\bar{n}+1,p}_B$ has weak derivatives of order $[ \frac{\bar{n}+1}{2i+1}]$ in
any direction of the increments $[i]$, $i\le [\frac{\bar n}{2}]$, it is not difficult to check,
similarly to \eqref{e11b} that $$\|F_{11}(\cdot,\z)\|_p\lesssim
|\x|_B^{\bar{n}+1}\|u\|_{W^{n+1,p}_B}.$$ On the other hand, by \eqref{ae16}, we have
\begin{align}
 F_{12}(z,\z)
 =\sum_{\bbB\le \bar n}\frac{\x^{\b}}{\b!}\left(\p^{\b}u(t,\bar x)-\p^{\b}u(t, x)\right).
\end{align}
Then, taking the $L^{p}$ norm in $dz$, we have
\begin{align}
 \|F_{12}(\cdot,\z)\|_p&\lesssim \sum_{\bbB\le \bar n}\frac{{|\x|_{B}^{\bbB}}}{\b!} \|\p^{\b}u(\cdot
 \circ 
 {\x^{[[\frac{\bar  n+1}{2}]:r]}})-\p^{\b}u\|_p
\end{align}
Now we use the inductive hypothesis on $\p^\b u\in W^{\bar{n}-\bbB+1,p}_B$ for $|\b|_{B}\ge 1$,
and finally get
  $$\|F_{12}(\cdot,\z)\|_p\lesssim \sum_{\bbB\le \bar n}\frac{{|\x|_{B}^{\bbB}}}{\b!}|\x|_B^{\bar{n}-\bbB+1}\|\p^{\b}u\|_{W^{\bar{n}-\bbB+1,p}_B}\lesssim
  {|\x|_{B}^{\bar{n}+1}}\|u\|_{W^{n+1,p}_B}.$$
\endproof

\begin{remark}
By Theorem \ref{ap1}, for any $i\in \{\bar d_{j-1}+1,\dots, \bar d_j\}$ with $2j+1>n$ we have in
particular
 $$
 \|u(e^{\d \p_{x_{i}}}\cdot)-u\|_p\lesssim |\d|^{\frac{n+1}{2j+1}}\|u\|_{W^{n+1,p}_B},\qquad u\in
 W^{n+1,p}_B.$$
Then, using Fubini's Theorem it is straightforward to check that, for any $\e>0$
  $$[u]_{\p_{x_i},\frac{n+1}{2j+1}-\e,p}{\lesssim {\e}^{-1}}\|u\|_{W^{n+1,p}_B},\qquad u\in W^{n+1,p}_B.$$
Together with Corollary \ref{P02} this gives the expected regularity in any spatial direction,
which is not prescribed a priori by the definition of the spaces. Also, by Corollary \ref{P02} we
can further infer that, for any $k\in\N_{0}$, $\b\in \N_0^N$ with $2j+1>n-2k-\bbB\ge 0$ we have
  $$[Y^k\p^{\b}u]_{\p_{x_i},\frac{n+1-2k-\bbB}{2j+1}-\e,p}{\lesssim {\e}^{-1}}\|u\|_{W^{n+1,p}_B},\qquad
  u\in W^{n+1,p}_B.$$
\end{remark}

\section{Approximation and interpolation}\label{Interpolation}
\subsection{Approximation in $W^{n,p}_B$}
Let $\phi$ be a test function supported on $\|z\|_B\le 1$ with unitary integral. Following
\cite{P22}, we define the $n$-th order approximation for $u\in W^{n,p}_{B}$ as
\begin{align}\label{approx}
  u_{n,\eps}(z):=
  \int_{\R^{N+1}}T_{n-1}u(\z,z)\phi\left(D_{\eps^{-1}}(\z^{-1}\circ
  z)\right)\frac{d\z}{\eps^{\Q}},\qquad \e>0,
\end{align}
where $T_nu(\z,z)$ is the $B$-Taylor polynomial in \eqref{eq:def_Tayolor_n} and $\Q$ the
homogeneous dimension of $\R^{N+1}$. Notice that
\begin{equation}\label{mollifier}
 \int_{\R^{N+1}}\phi\left(D_{\eps^{-1}}(\z^{-1}\circ z)\right)\frac{d\z}{\eps^{\Q}}=\int_{\|\z\|_B\le 1}\phi(\z)d\z=1.
\end{equation}
We also recall the useful Lemma 3.2 from \cite{P22} which still holds for functions in
$W^{n,p}_{B}$: indeed, its proof relies only on basic algebraic rules of derivation, namely the
Leibniz formula and the chain rule for compositions with smooth functions.
\begin{lemma}\label{LEM}
For any $u\in W^{n,p}_B$ and ${z,\z}\in \R^{N+1}$ we have
\begin{align}
 \p_{x_i}T_nu(\z,z)&=T_{n-1}(\p_{i}u)(\z,z), \qquad n\ge 1, \ i=1,\dots,d, \label{Taylor1}\\
 Y_zT_nu(\z,z)&=T_{n-2}(Yu)(\z,z), \qquad n\ge 2. \label{Taylor2}
\end{align}
\end{lemma}
\endproof
\begin{theorem}[\bf Approximation]\label{Approximation}
Let $n,m \in \N$ with $n<m$. {There exist constants $c_1=c(n,p,B)$ and $c_2=c(n,m,p,B)$} such that
for any $u\in W^{n,p}_B$ and $0<\eps \le 1$ we have
\begin{align}\label{Appr1}
 \|u-u_{n,\eps}\|_{p}& \le c_1 \eps^{n}\|u\|_{W^{n,p}_B}, \\ \|u_{n,\eps}\|_{W^{m,p}_B}&\le c_2
 \eps^{n-m}\|u\|_{W^{n,p}_B}.\label{Appr2}
\end{align}
\end{theorem}
\proof We denote by $\DB^l$ any weak derivative of intrinsic order $l$, that is
$\DB^l=Y^k\p_{x}^\b$ with $2k+\bbB=l$, and let
\begin{align}
 I^{(n,l)}_\e u(z):=\int_{\R^{N+1}}(T_{n-1}u(\z,z)-u(z)){(\DB^l\phi)}\left(D_{\e^{-1}}(\z^{-1}\circ z)\right)\frac{d\z}{\e^{\Q}}.
\end{align}
We prove the following preliminary estimates {for $u\in W^{n,p}_B\cap C^{\infty}_{0}$}:
\begin{align}
 \|I^{(n,l)}_\e u\|_{p}&\lesssim \e^{n-l}\|u\|_{W^{n,p}_B},\label{Appr3}\\
 [u_{0,\e}]_{Y,\frac 12,p}&\lesssim \e^{-1}\|u\|_p,\label{Appr4}\\ [I^{(n,l)}_\e u]_{Y,\frac 12,p}&\lesssim
 \e^{n-l-1}\|u\|_{W^{n,p}_B}.\label{Appr5}
\end{align}
First observe that by \eqref{invariance_law}, \eqref{homogeneity} and the change of variable
$\y=D_{\e^{-1}}(\z^{-1}\circ z)$, we have
\begin{align}
 I^{(n,l)}_\e
 u(z):&=\e^{-l}\int_{\R^{N+1}}\left(T_{n-1}u(\z,z)-u(z)\right)(\DB^l\phi)\left(D_{\e^{-1}}(\z^{-1}\circ
 z)\right)\frac{d\z}{\e^{\Q}}\\ &=\e^{-l}\int_{\|\y\|_B\le 1}\left(T_{n-1}u(z\circ (D_\e
 \y)^{-1},z)-u(z)\right)\DB^l\phi(\y){d\y}.
\end{align}
Then, by Minkowski integral inequality and Theorem \ref{ap1}, we get
\begin{align}
 \|I^{(n,l)}_\e u\|^p_{p}&\le \e^{-lp}\int_{\|\z\|_B\le 1}\|T_{n-1}u(\cdot\circ (D_\e
 \z)^{-1},\cdot)-u\|^p_p {|\DB^l\phi(\z)|}{d\z}\le
\intertext{(since $\|D_\e \z\|_B=\e\|\z\|_B$)}
  &\le \e^{-lp}\int_{\|\z\|_B\le 1}\e^{np}\|\z\|_B^{np}\|u\|^p_{W^{n,p}_B} {|\DB^l\phi(\z)|}{d\z}\lesssim
 \e^{(n-l)p}\|u\|^p_{W^{n,p}_B}
\end{align}
which proves \eqref{Appr3}. By a similar argument, we have
\begin{align}
 [u_{0,\e}]^p_{Y,\frac 12,p}&=\int_{\R^{N+1}}\int_{|h|\le 1}\left|\int_{\R^{N+1}}u(\z)
 \left(\phi\left(D_{\e^{-1}}(\z^{-1}\circ e^{h Y}z)\right)-\phi\left(D_{\e^{-1}}(\z^{-1}\circ
 z)\right)\right)\frac{d\z}{\e^{\Q}}\right|^p\frac{dh}{|h|^{1+\frac{p}{2}}}dz\\ &=\int_{|h|\le
 1}\int_{\R^{N+1}}\left|\int_{\R^{N+1}}u(z\circ (D_\e \z)^{-1})\left(\phi(e^{\frac{h}{\e^2}
 Y}\z)-\phi(\z)\right)d\z\right|^p dz\frac{dh}{|h|^{1+\frac{p}{2}}}\le
\intertext{(by Minkowski inequality)}
 &\le {\|u\|_{p}^{p}}\left({\int_{|h|> \e^2}+\int_{|h|\le \e^2}}\right){\left(\int_{\R^{N+1}}\left|\phi(e^{\frac{h}{\e^2}
 Y}\z)-\phi(\z)\right|d\z\right)^p}\frac{dh}{|h|^{1+\frac{p}{2}}} =:I_1+I_2
\end{align}
{By the triangular and H\"older inequalities we have
\begin{align}
I_1\le \|u\|^p_p\int_{|h|>\e^2}\left(\int_{\|\z\|_B\le
1}2\phi(\z)d\z\right)^p\frac{dh}{|h|^{1+\frac p2}}\lesssim
\|u\|^p_p\|\phi\|_p^p\int_{|h|>\e^2}\frac{dh}{|h|^{1+\frac{p}{2}}}\lesssim
\e^{-p}\|u\|^p_p\|\phi\|_p^p.
\end{align}
Next, noting that $e^{\frac{h}{\e^2}Y}\z=\z\circ(\frac{h}{\e^2},0)$, by \eqref{normcontrol} we
have $$\|e^{\frac{h}{\e^2}Y}\z\|_B\le m\left(\|\z\|_B+\frac{\sqrt{|h|}}{\e}\right)\le 2m $$ in the
integration set of $I_2$, and therefore, again by  H\"older inequality, we have
\begin{align}
 I_2\lesssim \|u\|^p_p\int_{|h|\le
 \e^2}\int_{\R^{N+1}}\left|\phi(e^{\frac{h}{\e^2}Y}\z)-\phi(\z)\right|^p d\z \frac{dh}{|h|^{1+\frac
 p2}}\le \e^{-p}\|u\|^p_p[\phi]^p_{Y,\frac 12,p}
\end{align}
where the last inequality easily follows by a change of variables, and this proves \eqref{Appr4}.}

Lastly, \eqref{Appr5} requires more attention. We have
 $$[I^{(n,l)}_\e u]^p_{Y,\frac 12,p}\le S_1+S_2$$
where
\begin{align}
 S_1&=\int_{\R^{N+1}}\int_{|h|\le
 1}\bigg|\int_{\R^{N+1}}\left(T_{n-1}u(\z,e^{hY}z)-u(e^{hY}z)\right)\times\\ &\quad \qquad \times
 \left((\DB^l\phi)\left(D_{\e^{-1}}(\z^{-1}\circ
 e^{hY}z)\right)-(\DB^l\phi)\left(D_{\e^{-1}}(\z^{-1}\circ
 z)\right)\right)\frac{d\z}{\e^{\Q}}\bigg|^p\frac{dh}{|h|^{1+\frac{p}{2}}}dz,\\
 S_2&= \int_{\R^{N+1}}\int_{|h|\le
 1}\Big|\int_{\R^{N+1}}J_{n}(\z,z) \DB^l\phi\left(D_{\e^{-1}}(\z^{-1}\circ
 z)\right)\frac{d\z}{\e^{\Q}}\Big|^p\frac{dh}{|h|^{1+\frac{p}{2}}}dz
\intertext{with}
 J_{n}(\z,z)&=T_{n-1}u(\z,e^{hY}z)-u(e^{hY}z)-\left(T_{n-1}u(\z,z)-u(z)\right).
\end{align}
The term $S_1$ can be controlled as $[{u}_{0,\e}]_{Y,\frac 12,p}$: {for simplicity here we assume
that we can control the support of the increment of $\DB^l\phi$ independently of $\e$ on the whole
integration set of $h$ to exploit a H\"older inequality (otherwise we can just split the
integration set and proceed as for the terms $I_1$ and $I_2$ of $[{u}_{0,\e}]_{Y,\frac 12,p}$),
then we have}
\begin{align}
 S_1&=\e^{-lp}\int_{|h|\le 1}\int_{\R^{N+1}}\Big|\int_{\R^{N+1}}\left(T_{n-1}u(e^{hY}z\circ (D_\e
 \z)^{-1},e^{hY}z)-u(e^{hY}z)\right)\times\\ &\quad \qquad \times \left(\DB^l\phi(e^{\frac{h}{\e^2}
 Y}\z)-\DB^l\phi(\z)\right)d\z\Big|^pdz\frac{dh}{|h|^{1+\frac{p}{2}}}
 \\&{\le}
 \e^{-lp}\int_{|h|\le 1}\int_{\R^{N+1}}\|T_{n-1}u(\cdot\circ (D_\e
 \z)^{-1},\cdot)-{u}\|^p_p|\DB^l\phi(e^{\frac{h}{\e^2}
 Y}\z)-\DB^l\phi(\z)|^pd\z\frac{dh}{|h|^{1+\frac{p}{2}}}\lesssim
\intertext{(by Theorem \ref{ap1})}
  &{\lesssim}
 \e^{(n-l)p}\|u\|^p_{W^{n,p}_B}\e^{-p}[\DB^l\phi]^p_{Y,\frac 12,p}\lesssim
 e^{(n-l-1)p}\|u\|^p_{W^{n,p}_B}.
\end{align}
On the other hand, we have $S_{2}=S_{21}+S_{22}$ where
\begin{align}
 S_{21}&:= \int_{\R^{N+1}}\int_{|h|\le\e^2}\Big|\int_{\R^{N+1}}J_{n}(\z,z)\e^{-l}(\DB^l\phi)\left(D_{\e^{-1}}(\z^{-1}\circ
  z)\right)\frac{d\z}{\e^{\Q}}\Big|^p\frac{dh}{|h|^{1+\frac{p}{2}}}dz\le
\intertext{(reasoning as above)}
 &\le \e^{-lp} \int_{\R^{N+1}}\int_{|h|>\e^2}2\|T_{n-1}u(\cdot\circ (D_\e
 \z)^{-1},\cdot)-u\|^p_p|\DB^l\phi(\z)|^p \frac{dh}{|h|^{1+\frac{p}{2}}}d\z \\ &\lesssim
 \e^{(n-l)p}\|u\|^p_{W^{n,p}_B}\|\DB^l\phi\|^p_p\int_{|h|>\e^2}\frac{dh}{|h|^{1+\frac{p}{2}}}\lesssim
 \e^{(n-l-1)p}\|u\|^p_{W^{n,p}_B}
\intertext{and}
 S_{22}&:= \int_{\R^{N+1}}\int_{\e^2<|h|\le
 1}\Big|\int_{\R^{N+1}}J_{n}(\z,z)\e^{-l}(\DB^l\phi)\left(D_{\e^{-1}}(\z^{-1}\circ
  z)\right)\frac{d\z}{\e^{\Q}}\Big|^p\frac{dh}{|h|^{1+\frac{p}{2}}}dz.
\end{align}
To estimate $S_{22}$, assume for a moment that $n\ge 3$: then by Theorem \ref{ap1} we have
\begin{align}
 S_{22}&\le \e^{-lp} \int_{|h|\le\e^2}\Big(\int_{\|\z\|_B\le
 1}\Big(\int_{\R^{N+1}}|h|^p\Big[\int_{0}^1|T_{n-3}Yu(z\circ (D_{\e}\z)^{-1},e^{\l
 hY}z)+\\ &\quad \quad -Yu(e^{\l hY}z)|d\l\Big]^pdz\Big)^{\frac
 1p}|\DB\phi(\z)|d\z\Big)^p\frac{dh}{|h|^{1+\frac{p}{2}}}\label{S21}.
\end{align}
Notice that
  $$z\circ (D_{\e}\z)^{-1}=e^{\l hY}z\circ (-\l h,0)\circ (D_{\e}\z)^{-1}=e^{\l hY}z\circ (e^{\l hY}D_{\e}\z)^{-1}.$$
Then, after the change of variables $\bar z=e^{\l hY}z$ and exchanging the order of integration,
the term inside the square brackets in \eqref{S21} is bounded by
\begin{align}
 \int_0^1\|T_{n-3}Yu(\cdot\circ (e^{\l hY}D_{\e}\z)^{-1},\cdot)-Yu\|^p_pd\l
 \lesssim 
 \int_0^1\|e^{\l hY}D_{\e}\z\|^{(n-2)p}_p\|Yu\|^p_{W^{n-2,p}_B}d\l.
\end{align}
By \eqref{homogeneity} and \eqref{normcontrol2}, recalling that $\l h\le \e^2$ and $\|\z\|_B\le 1$
in the current integration set, we have
 $$\|e^{\l hY}D_{\e}\z\|_B=\|D_{\e}(e^{\l \frac{h}{\e^2}Y}\z)\|_B=\e\|e^{\l \frac{h}{\e^2}Y}\z\|_B\lesssim \e (1+\|\z\|_B)\lesssim \e.$$
Therefore, substituting in \eqref{S21} we find
\begin{align}
 S_{22}&\le \e^{-lp}\|Yu\|^p_{W^{n-2,p}_B}\int_{|h|\le\e^2}\left(\int_{\|\z\|_B\le
 1}\e^{(n-2)}|\DB\phi(\z)|d\z\right)^p\frac{dh}{|h|^{1-\frac{p}{2}}}\\ &\lesssim
 \e^{(n-2-l)p}\|Yu\|^p_{W^{n-2,p}_B}\|\DB^l\phi\|_p^p\int_{|h|\le\e^2}\frac{dh}{|h|^{1-\frac{p}{2}}}\lesssim
 \e^{(n-1-l)p}\|u\|^p_{W^{n,p}_B}.
\end{align}
The cases $n=1$ or $n=2$ are easier: it is easy to check that $T_0u(\z,e^{hY}z)=T_0u(\z,z)=u(\z)$
and $T_1u(\z,e^{hY}z)=T_1u(\z,z)$, therefore it suffices to use \eqref{e11b} and proceed as above.
Collecting the estimates for $S_1$, $S_{21}$, $S_{22}$ we get \eqref{Appr5}.

\medskip
We are ready to prove \eqref{Appr1} and \eqref{Appr2} {for $u\in W^{n,p}_B\cap C^{\infty}_0$, then
the general statement follows by density}. Clearly $$\|u-u_{n,\e}\|_p=\|I^{(n,0)}_\e u\|_p\lesssim
\e^{n}\|u\|_{W^{n,p}_B}$$ by \eqref{Appr3}. On the other hand, by \eqref{eqa5}, with some slight
abuse of notation, we have
\begin{equation}\label{Appr6}
 |u_{n,\e}|_{m,p,B}\lesssim \|\DB^m u_{n,\e}\|_p+[\DB^{m-1} u_{n,\e}]_{Y,\frac 12,p}.
\end{equation}
{Since $\DB^i_zT_nu(\z,z)=0$ for any $i>n$ we have
\begin{align}
 \DB^mu_{n,\e}(z)=\sum_{i=0}^{{n}} \int_{\R^{N+1}}\DB^i_zT_{n-1}u(\z,z)\DB^{m-i}_z\phi(D_{\e^{-1}}(\z^{-1}\circ z))\frac{d\z}{\e^{\Q}},
\end{align}
meaning that $\DB^m$, $\DB^{m-i}$, $\DB^i$ may stand for any intrinsic derivative of order $m,
m-i, i$. Then, using Lemma \ref{LEM} and that $$\int_{\R^{N+1}}\DB^i\phi(z)dz=0, \qquad i>1, $$ we
can write
\begin{align}
  \DB^mu_{n,\e}(z)=\sum_{i=0}^n\int_{\R^{N+1}}\left(T_{n-1-i}\DB^iu(\z,z)-\DB^{i}u(z)\right)\DB^{m-i}_z\phi(D_{\e^{-1}}(\z^{-1}\circ
z))\frac{d\z}{\e^{\Q}},
\end{align}
and thus
\begin{equation}\label{Appr7}
\|\DB^m u_{n,\e}\|_p\lesssim \sum_{i=0}^n\|I^{(n-i,m-i)}_{\e}\DB^i u\|_p\lesssim \sum_{i=0}^n
\e^{n-m}\|\DB^i u\|_{W^{n-i,p}_B}\lesssim \e^{n-m} \|u\|_{W^{n,p}_B}.
\end{equation}}
It only remains to estimate $[\DB^{m-1} u_{n,\e}]_{Y,\frac 12,p}$: as before, if $m>n+1$, the test
function is affected by at least one derivative for any non-null term of $\DB^{m-1} u_{n,\e}$.
Therefore we have
\begin{align}
 [\DB^{m-1} u_{n,\e}]_{Y,\frac 12,p}\lesssim\begin{cases} \sum\limits_{i=0}^n[I^{(n-1,m-1-i)}_\e \DB^i
 u]_{Y,\frac 12,p}, & m>n+1\\ \sum\limits_{i=0}^{n-1}[I^{(n-1,m-1-i)}_\e \DB^i u]_{Y,\frac
 12,p}+[(\DB^nu)_{0,\e}]_{Y,\frac 12,p},  & m=n+1.
\end{cases}
\end{align}
By \eqref{Appr4} and \eqref{Appr5}, we directly derive
\begin{equation}\label{Appr8}
[\DB^{m-1} u_{n,\e}]_{Y,\frac 12,p}\lesssim \e^{n-m} \|u\|_{W^{n,p}_B},
\end{equation}
and recalling \eqref{Appr6}, \eqref{Appr7} we finally get \eqref{Appr2}.
\endproof

\subsection{Interpolation on the degree of smoothness}
{In this section, we establish an interpolation result. The primary definitions and key results
pertaining to interpolation theory are succinctly summarized in Appendix \ref{interpolation}.}
\begin{theorem}[\bf Interpolation]\label{t3}
For $1\le n\le m$ and {$1\le p\le \infty$} we have
\begin{equation}\label{ae22}
  (L^p,W^{m,p}_B)_{\frac{n}{m},1}\subseteq W^{n,p}_B\subseteq (L^p,W^{m,p}_B)_{\frac{n}{m},\infty}.
\end{equation}
\end{theorem}
\proof The first embedding in \eqref{ae22} is a direct consequence of Proposition \ref{p1}.
Indeed, from \eqref{interp_1} we deduce
\begin{equation}\label{interp_3}
  \|u\|_{W^{n,p}_B}\lesssim \eps \|u\|_{W^{m,p}_B}+\eps^{-\frac{n}{m-n}}\|u\|_p,\qquad \e>0.
\end{equation}
In particular, taking the optimal
$\eps=\left({\|u\|_p}/{\|u\|_{W^{m,p}_B}}\right)^{\frac{m-n}{m}}$ we get
\begin{equation}\label{interp_4}
 \|u\|_{W^{n,p}_B}\lesssim \|u\|_{W^{m,p}_B}^{\frac{m-n}{m}}\|u\|_p^{\frac nm}
\end{equation}
and, by Proposition \ref{p2}, estimate \eqref{interp_4} is equivalent to the embedding
$(L^p,W^{m,p}_B)_{\frac{n}{m},1}\subseteq W^{n,p}_B$.

The second embedding in \eqref{ae22} is a direct consequence of Theorem \ref{Approximation}.
Indeed, by \eqref{Appr1} and \eqref{Appr2}, for any $\e>0$ we have $$t^{-\frac
nm}K(t,u;L^p,W^{m,p}_B)\lesssim t^{-\frac nm}(\|u-u_{n,\e}\|_p+t\|u_{n,\e}\|_{W^{m,p}_B})\lesssim
t^{-\frac nm}(\e^n+t\e^{n-m})\|u\|_{W^{n,p}_B}.$$ Therefore, taking $\e=t^{\frac 1m}$ we get
$$t^{-\frac nm}K(t,u;L^p,W^{m,p}_B)\lesssim \|u\|_{W^{n,p}_B}.$$
\endproof


\section{Proof of Theorem \ref{Embedd_1}}\label{proofmain}
The proof of Theorem \ref{Embedd_1} is based on some basic {result} from interpolation theory,
which we briefly recall in Appendix \ref{interpolation} for reader's convenience. As a first step
we derive a ``crude'' embedding which will serve as a starting point to derive the general result,
through the characterization of Lorentz spaces of Lemma \ref{lem_Lorentz}.
\begin{lemma}\label{lem1}
For $p\in [1,\infty)$, let $u\in W^{1,p}_B$. There exists $r={r(p,B)}>p$, such that
\begin{equation}\label{emb3}
\|u\|_q\lesssim \|u\|_p^{1-\theta}|u|_{1,p,B}^{\theta}, \qquad q\in [p,r), \; \theta=\Q
\left(\frac{1}{p}-\frac{1}{q}\right).
\end{equation}
\end{lemma}
\proof As in the proof of {Theorem \ref{t3}} we start from the decomposition
$u=u_{\eps}+(u-u_{\eps})$ with $u_{\e}= u_{\e,n}$ as in \eqref{approx} with $n=1$ and $\e>0$.
By \eqref{Appr1}
  $$\|u-u_\eps\|_p\lesssim \eps \|u\|_{W^{1,p}_B}.$$
On the other hand, by Young's inequality we have
  $$\|u_{\eps}\|_{\infty}\le 
  \eps^{-\Q }\|u\|_p\|\phi(D_{\eps^{-1}})\|_{{\frac{p}{p-1}}}\lesssim
  \eps^{-\frac{\Q}{p}}\|u\|_p,\qquad \e>0.$$
Therefore, for $K$ as in \eqref{ae26}, we have
  $$K(t,u,L^p,L^{\infty})\lesssim (\eps+t \eps^{-\frac{\Q}{p}})\|u\|_{W^{1,p}_B},\qquad \e,t>0.$$
In particular, for $\eps=t^{\frac{p}{\Q+p}}$ we get
  $${K(t,u,L^p,L^{\infty})\lesssim t^{\frac{p}{\Q+p}}\|u\|_{W^{1,p}_B},\qquad t>0,}$$
that is
  $$W^{1,p}_B\subseteq (L^p,L^{\infty})_{\frac{p}{\Q+p},\infty}={L^{r}_{w}},\qquad r=\frac{p(\Q+p)}{\Q}>p,$$
by \eqref{Lorentz_K}. {By \eqref{ae24}}, this yields in particular that
$W_B^{1,p}\subseteq L^q$, for any $q\in [p,r)$ and
\begin{equation}\label{emb1}
  \|u\|_q\lesssim \|u\|_p+|u|_{1,p,B}, \qquad u\in W^{1,p}_B.
\end{equation}
Finally, by the usual scaling argument, applying \eqref{emb1} to $u(D_{\eps^{-1}}(\cdot ))$, and
using {\eqref{scaling}} we get
\begin{equation}\label{emb2}
 \|u\|_q\lesssim
  \eps^{\Q(\frac{1}{p}-\frac{1}{q})}\|u\|_p+\eps^{\Q(\frac{1}{p}-\frac{1}{q})-1}|u|_{1,p,B},
 \qquad u\in W^{1,p}_B, \; \eps>0,
\end{equation}
and this directly yields \eqref{emb3}, optimizing on $\eps$.
\endproof
{\begin{remark}
For $1\le p<\Q$ we must have $q\le p^{\ast}=\frac{\Q p}{\Q-p}$ in \eqref{emb2} or we would get a
contradiction by letting $\e$ tend to $0$. This means that the critical exponent $p^{\ast}$ is
optimal for the embedding \eqref{Embedd_1e2}.
\end{remark}
}


Next we use an ingenious approach, devised by Tartar \cite{MR1662313}, which consists of applying
\eqref{emb3} to a suitable non-linear transformation of $u$. Precisely, we consider $\phi_{k}(u)$
where $(\phi_k)_{k\in\Z}$ is an appropriate sequence of functions involving the levels
$a_k:=u^{\ast}(e^{k})$ of Lemma \ref{lem_Lorentz}: for $v\in\R$ and $k\in\Z$ we set
\begin{align}
 \phi_k(v)=\begin{cases}
 0& \text{ if }\ |v|\le a_{k+1},\\
 |v|-a_{k+1}& \text{ if }\ a_{k+1}\le |v|\le a_{k},\\
 a_{k}-a_{k+1} & \text{ if }\ |v|\ge a_{k}.
\end{cases}
\end{align}
We have the following crucial
\begin{lemma}
For $p\in [1,\infty)$, let $u\in W^{1,p}_B$. There exists a positive constant $c={c(p,B)}$ such
that
 \begin{equation}\label{emb5}
 e^{\frac{k}{p^{\ast}}}(a_k-a_{k+1})\le c 
 |\phi_k(u)|_{1,p,B},\qquad k\in\Z,
 \end{equation}
where $p^\ast$ is the critical exponent in \eqref{Embedd_1e2}.
\end{lemma}
\proof Notice that
\begin{equation}
  (a_{k}-a_{k+1})\caratt_{(|u|\ge a_{k})}\le \phi_{k}(u)\le (a_{k}-a_{k+1})\caratt_{(|u|\ge
  a_{k+1})},
\end{equation}
where $\caratt_{A}$ denotes the indicator function of the set $A$. Hence, for any $q\ge 1$ we have
\begin{align}
  \text{\rm Leb}(|u|\ge a_{k})^{\frac{1}{q}}(a_{k}-a_{k+1})\le \|\phi_{k}(u)\|_{q}\le(a_{k}-a_{k+1})\text{\rm Leb}(|u|\ge
  a_{k+1})^{\frac{1}{q}}
\end{align}
{where $\text{\rm Leb}(\cdot)$ represents the Lebesgue measure.} By \eqref{ae27}, which follows by
construction of $(a_{k})_{k\in\Z}$ (also recall definition \eqref{ae28} of distribution function),
we get
\begin{equation}\label{ae29}
  e^{\frac{k-1}{q}}(a_{k}-a_{k+1})\le \|\phi_{k}(u)\|_{q}\le e^{\frac{k+1}{q}}(a_{k}-a_{k+1})
\end{equation}
From the first inequality in \eqref{ae29} and \eqref{emb3} applied to $\phi_{k}(u)\in W^{1,p}_B$
with $q,\theta$ as in Lemma \ref{lem1}, we infer
\begin{align}
  e^{\frac{k-1}{q}}(a_{k}-a_{k+1})&\lesssim
  \|\phi_{k}(u)\|_{p}^{1-\theta}|\phi_{k}(u)|^{\theta}_{1,p,B}\lesssim
\intertext{(by the second inequality in \eqref{ae29})}
  &\lesssim
  e^{\frac{(k+1)(1-\theta)}{p}}(a_{k}-a_{k+1})^{1-\theta}|\phi_{k}(u)|^{\theta}_{1,p,B}.
\end{align}
Equivalently, we have
  $$e^{k\left(\frac{1}{q}-\frac{1-\theta}{p}\right)}(a_{k}-a_{k+1})^{\theta}\lesssim
  e^{\frac{k-1}{q}-\frac{(k+1)(1-\theta)}{p}}(a_{k}-a_{k+1})^{\theta}\lesssim |\phi_{k}(u)|^{\theta}_{1,p,B}$$
and this concludes the proof since $\frac{1}{q}-\frac{1-\theta}{p}=\frac{\theta}{p^{\ast}}$.
\endproof

\proof[Proof of Theorem \ref{Embedd_1}] Using that $|\phi'_{k}(v)|=1$ for $a_{k+1}< |v|< a_{k}$
and $\phi'_{k}(v)=0$ elsewhere, it is not difficult to prove that
\begin{equation}\label{emb4}
  |u|_{1,p,B}<\infty\quad \text{ if and only if }\quad |\phi_k(u)|_{1,p,B}\in \ell^p(\Z).
\end{equation}
Thus, combining \eqref{emb5} with \eqref{emb4}, we deduce that
$e^{\frac{k}{p^{\ast}}}(a_k-a_{k+1})\in \ell^p(\Z)$ for any  $u\in W^{1,p}_{B}$.

\medskip\noindent{[\it Case $1\le p<\Q$}] A direct application of Lemma \ref{lem_Lorentz} gives the improved
Sobolev embedding
\begin{equation}
  W^{1,p}_B\subseteq L^{p^{\ast},p}\subseteq L^{p^{\ast}}.
\end{equation}
In particular $W^{1,p}_B\subseteq L^{q}$ for any $q\in [p,p^\ast]$ by a standard application of
the Young inequality.


\medskip\noindent[{\it Case $p>\Q$}] We have $p^{\ast}<0$ and therefore for any $\bar{k}\le 0$ we have
\begin{align}
  a^{p}_{\bar{k}}&=\left(\sum_{k=-\infty}^{\bar{k}}(a_k-a_{k+1})\right)^{p}\\ &\le
  \sum_{k=-\infty}^{\bar{k}}(a_k-a_{k+1})^{p}e^{\frac{pk}{p^{\ast}}}
  \left(\sum_{k=-\infty}^{\bar{k}}e^{-\frac{pk}{(p-1)p^{\ast}}}\right)^{p-1}\lesssim
  \sum_{k\in\Z}(a_k-a_{k+1})^{p}e^{\frac{pk}{p^{\ast}}}<\infty
\end{align}
by \eqref{emb5}.
Being decreasing, $(a_{k})_{k\in\Z}$ is then a bounded sequence and this yields
$W^{1,p}_B\subseteq L^{\infty}$, that is $\|u\|_{\infty}\lesssim \|u\|_p+|u|_{1,p,B}$: by the
usual scaling argument we have
\begin{equation}\label{emb7}
  \|u\|_{\infty}\lesssim \|u\|_p^{1-\frac{\Q}{p}}|u|_{1,p,B}^{\frac{\Q}{p}}.
\end{equation}
Applying \eqref{emb7} to $(u(e^{h \p_{x_i}}\cdot)-u)$, for $i=1,\dots,d$, using \eqref{e10} and
noticing that $|{u(e^{h \p_{x_i}}\cdot)}-u|_{1,p,B}\le 2|u|_{1,p,B}$, we get
\begin{equation}
 \|{u(e^{h \p_{x_i}}\cdot)}-u\|_{\infty}\lesssim
 |h|^{1-\frac{\Q}{p}}\|\p_{x_i}u\|_p^{1-\frac{\Q}{p}}|u|_{1,p,B}^{\frac{\Q}{p}} \lesssim
 |h|^{1-\frac{\Q}{p}}|u|_{1,p,B}, \qquad i=1,\dots,d.
\end{equation}
Analogously, applying \eqref{emb7} to $(u(e^{h Y}\cdot)-u)$ and using \eqref{e11} we get
\begin{equation}
 \|{u(e^{h Y}\cdot)}-u\|_{\infty}\lesssim |h|^{\frac{1}{2}(1-\frac{\Q}{p})}[u]^{1-\frac{\Q}{p}}_{Y,\frac
 12,p}|u|_{1,p,B}^{\frac{\Q}{p}} \lesssim |h|^{\frac{1}{2}(1-\frac{\Q}{p})}|u|_{1,p,B},
\end{equation}
which proves the Morrey embedding \eqref{Embedd_1e1}.

\medskip\noindent[{\it Case $p=\Q$}] 
As in the case $p<\Q$, the embeddings \eqref{ae31} follow from Lemma \ref{lem_Lorentz}. To get
estimate \eqref{ae30}, it suffices to repeat the argument used by Tartar in \cite{MR2328004},
Chapter 30: more precisely, for $p=\Q$ we have $1/p^{\ast}=0$ so that
$a_{k}-a_{k+1}\in \ell^\Q(\Z)$ by \eqref{emb5}-\eqref{emb4}; applying H\"older's inequality 
we first prove that for every $\e>0$ there exists a constant $c=c(\e,u)>0$ such that
  $$a_{k}^{\frac{\Q}{\Q-1}}\le \e |k|+c,\qquad k\le 0.$$
On the set where $a_{k+1}\le |u|<a_{k}$, which has measure less than $e^{k+1}$ by \eqref{ae27}, we
have
  $$e^{\l|u|^{\frac{\Q}{\Q-1}}}\le e^{\l a_{k}^{\frac{\Q}{\Q-1}}}\le e^{\l(\e |k|+c)},\qquad k\le 0,$$
and by choosing $\e<\tfrac{1}{\l}$ we deduce that $e^{\l|u|^{\frac{\Q}{\Q-1}}}$ is integrable on any
set where $|u|\ge \d>0$.
\endproof


\section{Higher orders embeddings}\label{Embedding_section}
{Embeddings for higher order Sobolev spaces are classically derived by iteration from the $n=1$
case. In our setting, because of the qualitative difference between even and odd orders of
intrinsic spaces (cf. \eqref{eqa6}) resulting from the two-steps iterative definition, we need
some additional work at least when $n=2$, in order to control the Holder and Sobolev-Slobodeckij
quasi-norms involving the vector field $Y$, in the high and low summability cases respectively.
Our method here is based on the representation of a $W^{2,p}_B$ function by means of the
fundamental solution of a linear Kolmogorov operator with drift matrix $B$ (see
\eqref{representation} below), which only allows to derive the embeddings in the case $p>1$. } The
main result of this section is the following
\begin{theorem}[\bf $W^{n,p}_B$ embeddings]\label{Embedd_k}
Let $k\in \N$, $n\in \N_0$ and ${1<p<\infty}$.
\begin{itemize}
  \item[i)] If $kp<\Q$, then
\begin{equation}\label{Sobolev_n}
  W^{n+k,p}_B\subseteq W_B^{n,q}, \qquad p\le q\le p^{\ast}_{k},
  \qquad \frac{1}{p^{\ast}_{k}}=\frac{1}{p}-\frac{k}{\Q}.
\end{equation}
In particular, $W^{k,p}_B\subseteq L^{q}$ for $p\le q\le p^{\ast}_{k}$;

\item[ii)] if $kp>\Q>(k-1)p$, then
\begin{equation}\label{Morrey_n}
  W^{n+k,p}_B\subseteq C_B^{n,k-\frac{\Q}{p}}.
\end{equation}
\end{itemize}
\end{theorem}
%
%

\subsection{Fundamental solution of Fokker-Planck equations}
We recall some preliminary results about the fundamental solution of the Fokker-Planck operator
$\K$ in \eqref{operator}. 
\begin{proposition} H\"ormander's condition \eqref{horm} is equivalent to the fact {that} the matrix
 $$\cC_t=\int_0^te^{sB}{\begin{pmatrix} I_{d} & 0_{d\times (N-d)}\\ 0_{(N-d)\times d} &
 0_{(N-d)\times (N-d)} \end{pmatrix}} e^{s{B^\ast}}ds $$
is positive definite every $t>0$. In this case, the fundamental solution of $\K$ with pole at $0$
is given by
\begin{equation}\label{FS}
 \G(t,x)=
 \begin{cases}\frac{1}{\sqrt{(2\pi)^N\det{\cC_t}}}e^{-\frac{1}{2}\langle {\cC}^{-1}_t x,x\rangle}, & t>0 \\ 0 &t\le 0.
 \end{cases}
\end{equation}
The fundamental solution with pole at $\z$ is the left {translation} of $\G$ with respect to the
group law, that is $\G(\z^{-1}\circ \cdot)$.
\end{proposition}

We recall that $\Q$ denotes the homogeneous dimension of $\R^{N+1}$ defined in \eqref{homdim}.
Since
\begin{equation}\label{homogeneity2}
 \cC_{\l^2t}=\hD_\l \cC_{t}\hD_\l,
\end{equation}
the fundamental solution $\G$ is homogeneous of degree ${-\Q+2}$ with respect to
{$(D_\l)_{\l>0}$}. Similarly $\p_{x_i}\G$, $\p_{x_ix_j}\G$ for $i,j=1,\dots,d$ and $Y\G$ are
homogeneous of degrees {$-\Q+1$, $-\Q$ and $-\Q$} respectively.

Later we will exploit the global estimates of the following
\begin{lemma}\label{l2}
For every $z\in \R^{N+1}$ we have
\begin{equation}\label{FScontrol}
 \G(z)\lesssim \|z\|^{{-\Q+2}}_B,\qquad |Y\G(z)|\lesssim  \|z\|^{{-\Q}}_B.
\end{equation}
\end{lemma}
\proof A local version of \eqref{FScontrol} has been proven in \cite{DiFrancescoPolidoro} in the
general framework of non-homogeneous Fokker-Planck operators. In our setting we provide a more
direct proof. Let us only consider the second estimate in \eqref{FScontrol} for $t>0$: {since
$\K\G=0$ we have
 \begin{equation}\label{FScontrol0}
  |Y\G(z)|
  \le \frac{1}{2}{\sum_{i=1}^{d}|\p_{x_ix_i}}\G(z)|\le
  \frac{1}{2}\sum_{i=1}^{d}\left(\frac{1}{2}|(\cC^{-1}_tx)_i|^{2}+|(\cC^{-1}_t)_{ii}|\right)\G(z).
  \end{equation}
By \eqref{FS} an \eqref{homogeneity2} we have
\begin{equation}\label{FScontrol1}
\G(z)\lesssim t^{1-\frac{\Q}{2}}\exp\left(-\frac{1}{2}\langle
 \hD_{\frac{1}{\sqrt{t}}}\cC^{-1}_1 \hD_{\frac{1}{\sqrt{t}}}x,x\rangle\right)\lesssim t^{1-\frac{\Q}{2}}\exp\left(-\frac{1}{2\|\cC^{-1}_{1}\|}
 |\hD_{\frac{1}{\sqrt{t}}}x|^2\right);
 \end{equation}
On the other hand, for $1\le i\le d$, again by \eqref{homogeneity2} we have
\begin{align}
 |(\cC^{-1}_tx)_i|&=|(\hD_{\frac{1}{\sqrt{t}}}\cC^{-1}_1
 \hD_{\frac{1}{\sqrt{t}}}x)_i|
 \le
 \frac{\|\cC^{-1}_1\|}{\sqrt{t}}|\hD_{\frac{1}{\sqrt{t}}}x|,\label{FScontrol2}\\
 |(\cC^{-1}_t)_{ii}|&=|(\cC^{-1}_t\mathbf{e}_i)_i|\le \frac{\|\cC^{-1}_1\|}{{t}}. \label{FScontrol3}
\end{align}
Therefore, using that $\|z\|_B=t^{\frac 12}\|(1,D_{\frac{1}{\sqrt{t}}}x)\|_B= t^{\frac{1}{2}}(1+|\hD_{\frac{1}{\sqrt{t}}}x|_B)$, by \eqref{FScontrol0}, \eqref{FScontrol1} and \eqref{FScontrol2}-\eqref{FScontrol3} we finally get
\begin{align}
 \|z\|^{\Q}_B|Y\G(z)|&\lesssim t^{\frac{\Q}{2}}(1+|\hD_{\frac{1}{\sqrt{t}}}x|_B)^{\Q}
 \frac{1}{t^{\frac{\Q}{2}}}(1+|\hD_{\frac{1}{\sqrt{t}}}x|^2) \exp\left(-\frac{1}{2\|\cC^{-1}_{1}\|}|\hD_{\frac{1}{\sqrt{t}}}x|^2\right)
 \lesssim 1.
\end{align}}
The proof is completed.
\endproof

\subsection{Proof of Theorem \ref{Embedd_k}}
We first examine the regularity along the vector field $Y$ in both the high and low summability
case for $W^{2,p}_B$. We recall a theorem for convolution with a homogeneous kernel proved in
\cite{MR0290095}, Theorem 1, p.119.
\begin{theorem}\label{Stein}
For every $\a\in\, ]0,\Q[$  and $g\in L^p(\R^{N+1})$ with $p>1$, the function
 $$I_{\a}(g)(z)=\int_{\R^{N+1}}\frac{g(\z)}{\|\z^{-1}\circ z\|^{\Q-\a}_B}d\z,$$ is a.e. defined and
there exists $c=c(p,\a)>0$ such that
 $$\|I_\a(g)\|_q\le c\|g\|_p,\qquad \frac{1}{p}=\frac{1}{q}+\frac{\a}{\Q}.$$
\end{theorem}

\begin{proposition}\label{Y_Holder}
If $p>\Q$, then there exists $c=c(p,B)$ such that
\begin{equation}
 {\esssup_{z\in\R^{N+1}}}\, |u(e^{\delta Y}z)-u(z)|\le c |\delta|^{1-\frac{\Q}{2p}}\|u\|_{W^{2,p}_B}, \qquad u\in W^{2,p}_B, \quad \d\in \R.
\end{equation}
\end{proposition}
\proof By the definition of fundamental solution, for $u\in C^{\infty}_0$ we have the
representation
\begin{equation}\label{representation}
  u(z)=-\int_{\R^{N+1}}\G(\z^{-1}\circ z)\K u(\z)d\z.
\end{equation}
Since $\G$ is homogeneous of degree $-\Q+2$, by Theorem \ref{Stein} and a density argument we
deduce that \eqref{representation} holds a.e. for any $u\in W^{2,p}_B$ as well. Then we have
\begin{align}
 |u(e^{\delta Y}z)-u(z)|&\le \left(\int_{\|\z^{-1}\circ z\|_B\ge c\sqrt{|\delta|}}+
 \int_{\|\z^{-1}\circ z\|_B< c\sqrt{|\delta|}}\right)|\G(\z^{-1}\circ e^{\delta
 Y}z)-\G(\z^{-1}\circ z)| |\K u(\z)|d\z\le
\intertext{(for some $\bar{\d}$, $|\bar{\d}|\le |\delta|$, dependent on $z,\z$)}
 &\le \int_{\|\z^{-1}\circ z\|_B\ge c\sqrt{|\delta|}}|\delta|
 |Y\G (e^{\bar{\d}Y}(\z^{-1}\circ z))| |\K u(\z)|d\z\\ &\qquad +\int_{\|\z^{-1}\circ
 z\|_B<c\sqrt{|\delta|}} \left(\G(\z^{-1}\circ e^{\d Y}z)+\G(\z^{-1}\circ z)\right)|\K
 u(\z)|d\z\le S_1+S_{21}+S_{22}
\end{align}
by the uniform estimates \eqref{FScontrol}, where
\begin{align}
 S_1&= \int_{\|\z^{-1}\circ
 z\|_B\ge c\sqrt{|\delta|}}|\delta|\|e^{\bar{\d} Y}(\z^{-1}\circ
 z)\|_B^{-\Q} |\K u(\z)|d\z,\\
 S_{21}&=\int_{\|\z^{-1}\circ z\|_B< c\sqrt{|\delta|}}\|e^{\delta Y}(\z^{-1}\circ z)\|_B^{-\Q+2}|\K
 u(\z)|d\z,\\
 S_{22}&=\int_{\|\z^{-1}\circ z\|_B< c\sqrt{|\delta|}}\|\z^{-1}\circ z\|_B^{-\Q+2}|\K u(\z)|d\z.
\end{align}
Choosing $c$ as in \eqref{normcontrol2}, we have $\|\z^{-1}\circ z\|_B\lesssim \|\z^{-1}\circ
e^{\bar{\d} Y}z\|_B=\|e^{\bar{\d} Y}(\z^{-1}\circ z)\|_B$ on $\sqrt{|\delta|}\le c\|\z^{-1}\circ
z\|_B$. Then, as in \cite{MR1638177}, Lemma 2.9, we have
\begin{align}
 S_1&=|\delta|\sum_{k\ge 1}  \int_{c^k\sqrt{|\delta|}\le \|\z^{-1}\circ z\|_B\le
 c^{k+1}\sqrt{|\delta|}}\frac{|\K u(\z)|}{\|\z^{-1}\circ z\|_B^{\Q}}d\z\\
 &{\lesssim} |\delta|\sum_{k\ge 1}(c^k\sqrt{|\delta|})^{-\Q} \int_{\|\z^{-1}\circ
 z\|_B\le c^{k+1}\sqrt{|\delta|}}|\K u(\z)|{d\z}\le
\intertext{(by H\"older's inequality)}
 &\le |\delta|^{1-\frac{\Q}{2}}\sum_{k\ge 1} c^{-k\Q}\|\K u\|_p \text{\rm
 Leb}\left(\|\z^{-1}\circ z\|_B\le c^{k+1}\sqrt{|\delta|}\right)^{1-\frac 1p}\\
 &\lesssim
 |\delta|^{1-\frac{\Q}{2p}}\|u\|_{W^{2,p}_B}c^{\Q (1-\frac 1p)}\sum_{k\ge 1}{c^{-k\frac{\Q}{p}}},
\end{align}
where we used that $\text{\rm Leb}\left(\|\z^{-1}\circ z\|_B\le r\right)\le r^{\Q}\text{\rm
Leb}\left(\|z\|_B\le 1\right)$. Similarly, we have
\begin{align}
 S_{22}&=\sum_{k\ge 1} \int_{c^{-k}\sqrt{|\delta|}\le \|\z^{-1}\circ z\|_B\le
 c^{1-k}\sqrt{|\delta|}}\frac{|\K u(\z)|}{\|\z^{-1}\circ z\|_B^{\Q-2}}d\z\\ &\le \sum_{k\ge 1}
 \left(\frac{c^k}{\sqrt{|\delta|}}\right)^{\Q-2}\|\K u\|_p (c^{1-k}\sqrt{|\delta|})^{\Q
 (1-\frac 1p)} \lesssim |\delta|^{1-\frac{\Q}{2p}}\|u\|_{W^{2,p}_B}\sum_{k\ge
 1}c^{-k(2-\frac{\Q}{p})}.
\end{align}
For $S_{21}$ observe now that by \eqref{normcontrol2} on $\{\|\z^{-1}\circ z\|_B\le
c\sqrt{|\delta|}\}$ we have $\|{\z^{-1}}\circ e^{\delta Y}z\|_B\le
m(1+c)\sqrt{|\delta|}\equiv\bar{c}\sqrt{|\delta|}$, so that
  $$S_{21}\le \int_{\|\z^{-1}\circ e^{\delta Y}z\|_B\le \bar{c}\sqrt{|\delta|}}\|\z^{-1}\circ e^{\delta
  Y}z\|_B^{-\Q+2}{|\K u(\z)|}d\z,$$ and therefore it is analogous to $S_{22}$. The proof is complete.
\endproof

\begin{proposition}\label{Y_frac}
{If $p<\Q$, then there exists $c=c(p,B)$ such that
\begin{equation}
  [u]_{Y,\frac{1}{2},p^{\ast}}\le c \|u\|_{W^{2,p}_B}, \qquad u\in W^{2,p}_B.
\end{equation}}
\end{proposition}
\proof Recall that $$[u]^{p^{\ast}}_{Y,\frac{1}{2},p^{\ast}}=\int_{\R^{N+1}}\int_{|h|\le
1}\frac{|u(e^{h Y}z)-u(z)|^{p^{\ast}}}{|h|^{1+\frac{p^{\ast}}{2}}}dh dz.$$ Observe that, since
$|h|=\|z^{-1}\circ e^{hY}z\|_B=\|(e^{hY}z)^{-1}\circ z\|_B$, possibly exchanging variables by
$z'=e^{hY}z$ (whose Jacobian has determinant equal to $1$), we may assume that we are integrating
on a subset of
\begin{equation}\label{integration_set}
\{(z,h)\in \R^{N+1}\times [-1,1]\mid |u(z)|\ge |u(e^{h Y}z)|\}.
\end{equation}
Now, by representation \eqref{representation} and Minkowski inequality
\begin{align}
 [u]_{Y,\frac{1}{2},p^{\ast}}^{p^{\ast}}&=\int_{\R^{N+1}}\int_{|h|\le
 1}\Big|\int_{\R^{N+1}}\left(\G(\z^{-1}\circ e^{h Y}z)-\G(\z^{-1}\circ
 z)\right)\K u(\z)d\z\Big|^{p^{\ast}} \frac{dh}{|h|^{1+\frac{p^{\ast}}{2}}}dz\\
 &\le \int_{\R^{N+1}}\left(\int_{\R^{N+1}}\left(\int_{|h|\le
 1}|\G(\z^{-1}\circ e^{h Y}z)-\G(\z^{-1}\circ z)|^{p^{\ast}}
 \frac{dh}{|h|^{1+\frac{p^{\ast}}{2}}}\right)^{\frac{1}{p^{\ast}}} |\K
 u(\z)|d\z\right)^{p^{\ast}}dz.
\end{align}
For $c$ as in \eqref{normcontrol2}, we rewrite the inner integral of the last expression as
\begin{equation}
 \left(\int_{|h|\le c\|\z^{-1}\circ z\|^2_B}+\int_{c\|\z^{-1}\circ z\|^2_B\le |h|\le 1}\right)
 |\G(\z^{-1}\circ e^{h Y}z)-\G(\z^{-1}\circ
 z)|^{p^{\ast}}\frac{dh}{|h|^{1+\frac{p^{\ast}}{2}}}=:S_1+S_2.
\end{equation}
Now, for some $h'$ dependent on $h,z,\z$, with $|h'|\le |h|$, we have
\begin{align}
 S_1&=\int_{|h|\le c\|\z^{-1}\circ z\|^2_B}|Y\G(e^{h'Y}(\z^{-1}\circ
 z))|^{p^{\ast}}|h|^{p^{\ast}}\frac{dh}{|h|^{1+\frac{p^{\ast}}{2}}}\lesssim
\intertext{(by \eqref{FScontrol})}
 &\lesssim \int_{|h|\le
 c\|\z^{-1}\circ z\|^2_B}\|e^{h'Y}(\z^{-1}\circ z)\|_B^{-p^{\ast}\Q
 }\frac{dh}{|h|^{1-\frac{p^{\ast}}{2}}}\\
 &\lesssim \|\z^{-1}\circ z\|_B^{-p^{\ast}(\Q-1)}
\end{align}
using that $\|\z^{-1}\circ z\|_B\lesssim \|e^{h'Y}(\z^{-1}\circ
z)\|_B$ in the domain {of} the last integral. 

On the other hand, by \eqref{integration_set} and \eqref{FScontrol} we have
\begin{align}
 S_2&\le 2\int_{c\|\z^{-1}\circ z\|^2_B\le |h|}\G(\z^{-1}\circ
 z)^{p^{\ast}}\frac{dh}{|h|^{1+\frac{p^{\ast}}{2}}}\\ &\lesssim \|\z^{-1}\circ
 z\|^{-p^{\ast}(\Q-2)}\int_{c\|\z^{-1}\circ z\|^2_B\le |h|}\frac{dh}{|h|^{1+\frac{p^{\ast}}{2}}}
 \lesssim \|\z^{-1}\circ z\|_B^{-p^{\ast}(\Q-1)}.
\end{align}
Therefore we have
\begin{equation}
  [u]_{Y,\frac{1}{2},p^{\ast}}
  \lesssim \|I_1(\K u)|\|_{p^{\ast}}, 
\end{equation}
with $I_1$ as in {Theorem} \ref{Stein} and the thesis follows since $\K u\in L^{p}$ by assumption.
\endproof

We are now in position to prove Theorem \ref{Embedd_k}.

\proof[Proof of Theorem \ref{Embedd_k}] The embeddings of $W^{1,p}_{B}$ follow from Theorem
\ref{Embedd_1}. Regarding $W^{2,p}_{B}$, the statement of the theorem can be rewritten more
explicitly as follows:
\begin{itemize}
  \item[1)] if $p>\Q$ then
\begin{equation}\label{case3}
  W^{2,p}_{B}\subseteq C_{B}^{1,1-\frac{\Q}{p}};
\end{equation}

  \item[2)] if $\frac{\Q}{2}<p<\Q$ then:
\begin{itemize}
  \item considering $n=k=1$, we have
\begin{equation}\label{case2a}
  W^{2,p}_{B}\subseteq W^{1,q}_{B},\qquad p\le q\le p^{\ast},\qquad
  \frac{1}{p^{\ast}}=\frac{1}{p}-\frac{1}{\Q};
\end{equation}
  \item considering {$n=0$ and $k=2$}, we have
\begin{equation}\label{case2b}
  W^{2,p}_{B}\subseteq C_{B}^{0,2-\frac{\Q}{p}}.
\end{equation}
\end{itemize}
However, by \eqref{Embedd_1e1} we have that \eqref{case2a} implies \eqref{case2b} so that it
suffices to prove \eqref{case2a};

  \item[3)] if $p\le\frac{\Q}{2}$, which implies $n=0$ and $k=2$, we have
\begin{equation}\label{case1}
  W^{2,p}_{B}\subseteq L^{q},\qquad p\le q\le p^{\ast}_{2},\qquad
  \frac{1}{p^{\ast}_{2}}=\frac{1}{p}-\frac{2}{\Q};
\end{equation}
\end{itemize}

To prove \eqref{case3}, we notice that for $p>\Q$ we have $u, \p_{x_i}u\in W^{1,p}_B\subseteq
C_B^{0,1-\frac{\Q}{p}}$ for any $i=1,\dots,d$ by \eqref{Embedd_1e1} and
$[u]_{C_Y^{1-\frac{\Q}{2p}}}\lesssim \|u\|_{{W^{2,p}_B}}$ by Proposition \ref{Y_Holder}. 

To prove \eqref{case2a}, it suffices to observe that if $\frac{\Q}{2}<p<\Q$ then $u,\p_{x_i}u \in
W^{1,p}_B\subseteq L^{p^{\ast}}$ by \eqref{Embedd_1e2} and $[u]_{Y,\frac{1}{2},p^{\ast}}\lesssim
\|u\|_{W^{2,P}_B}$ by Proposition \ref{Y_frac}.


Finally, if $p<\frac{\Q}{2}$ then again $W^{2,p}_B\subseteq W^{1,p^{\ast}}_B$ with $p^{\ast}<\Q$
and by \eqref{Embedd_1e2} we get
  $$W^{2,p}_B\subseteq W^{1,p^{\ast}}_B\subseteq L^{\frac{p^{\ast}\Q }{\Q-p^{\ast}}}={L^{p^{\ast}_{2}}}.$$

The proof of higher order embeddings is analogous: by induction, it suffices to use iteratively
the previous arguments.
\endproof

\appendix\renewcommand{\theequation}{\thesection.\arabic{equation}}
\section{Interpolation}\label{interpolation}
We briefly recall some basic tool and notion from interpolation theory: for a comprehensive
presentation of the subject we refer, for instance, to \cite{MR3726909}, \cite{MR2328004} and
\cite{MR2424078}.

Given two real Banach spaces $Z_1, Z_2$, we write $Z_1=Z_2$ if $Z_1$ and $Z_2$ have the same
elements with equivalent norms; we write $Z_1\subseteq Z_2$ if $Z_1$ is continuously embedded in
$Z_2$. The pair $(Z_1, Z_2)$ is called an \textit{interpolation pair} if both $Z_1$ and $Z_2$ are
continuously embedded in some Hausdorff topological vector
space: in this case, the intersection $Z_1\cap Z_2$ and the sum 
$Z_1+Z_2$ endowed with the norms
  $$\|u\|_{Z_1\cap Z_2}:=\max\{\|u\|_{Z_1},\|u\|_{Z_2}\}, \qquad \|u\|_{Z_1+ Z_2}:=
  \inf_{u_1\in Z_1,\, u_2\in Z_2\atop u=u_{1}+u_{2}}(\|u_1\|_{Z_1}+\|u_2\|_{Z_2}),$$
are Banach spaces. For any $t>0$ and $u\in Z_{1}+Z_{2}$, we set
\begin{equation}\label{ae26}
 K(t,u)\equiv K(t,u;Z_1,Z_2):=
  \inf_{u_1\in Z_1,\, u_2\in Z_2\atop u=u_{1}+u_{2}}(\|u_1\|_{Z_1}+t\|u_2\|_{Z_2}).
\end{equation}
Any Banach space $E$ such that  
\begin{equation}\label{intermediate}
Z_1\cap Z_2\subseteq E\subseteq Z_1+Z_2,
\end{equation}
is called an \textit{intermediate space}. Among these, for $0<\theta<1$ and $1\le p\le \infty$, we
have the \textit{real interpolation space} $(Z_1,Z_2)_{\theta,p}$ consisting of $u\in Z_{1}+Z_{2}$
such that
\begin{equation}\label{k_interpolation}
  \|u\|_{\theta,p}:=\|t^{-\theta}K(t,u)\|_{L^{p}_{*}}<\infty
\end{equation}
where
  $L^{p}_{*}=L^{p}_{*}(\R_{>0})$
denotes the $L^{p}$ space with respect to the measure $\tfrac{dt}{t}$ on $\R_{>0}$ and $L^{\infty}_{*}:=L^{\infty}$. 
\begin{proposition}[\cite{MR3753604}, Prop. 1.20]\label{p2}
For an intermediate space $E$ the following conditions are equivalent:
\begin{itemize}
  \item[i)] $(Z_{1},Z_{2})_{\theta,1}\subseteq E$;
  \item[ii)] there exists a constant $c$ such that
\begin{equation}\label{ae21}
  \|u\|_{E}\le c\|u\|_{Z_{1}}^{1-\theta}\|u\|_{Z_{2}}^{\theta},\qquad u\in Z_{1}\cap Z_{2}.
\end{equation}
\end{itemize}
\end{proposition}
In the very particular case $Z_{2}\subseteq Z_{1}$ (for instance, if $Z_{1}$ is an $L^{p}$ space
and $Z_{2}$ is some intrinsic Sobolev space $W^{m,p}_{B}$), we have
\begin{equation}\label{ae20}
 Z_1\cap Z_2=Z_{2},\qquad Z_1+ Z_2=Z_{1},\qquad K(t,u)\le {\min\{\|u\|_{Z_1},t\|u\|_{Z_2}\}}.
\end{equation}
Then, since $t\mapsto K(t,u)$ is bounded by \eqref{ae20}, for $\|u\|_{\theta,p}$ in
\eqref{k_interpolation} to be finite, what really matters is only the behaviour of $K(t,u)$ near
$t=0$.

\subsection{Interpolation between $L^{p}$ spaces}\label{Lorentz_spaces}
The {\it distribution} of a measurable function  $u$ on $\R^{N}$ is defined as
\begin{align}\label{ae28}
  \m_{u}(\l)&:=\text{\rm Leb}(|u|>\l),\qquad \l\ge 0,
\intertext{while} \label{rearranging}
  u^{\ast}(t)&:=\inf\{\l\ge0\mid \m_{u}(\l)\le t \}, \qquad t\ge0,
\end{align}
is called the \textit{rearranging} of $u$. Distribution and rearranging are decreasing and right
continuous functions. Since
  $$u^{\ast}(t)>\l\qquad  \text{\rm if and only if}\qquad \m_{u}(\l)>t$$
we have
\begin{equation}\label{equimeasurability}
 \text{\rm Leb}(u^{\ast}>\l)=\text{\rm Leb}\left(0\le t<\m_{u}(\l)\right)=\text{\rm
 Leb}(|u|>\l)
\end{equation}
i.e. $u$ and $u^{\ast}$ are equimeasurable and consequently $\|u\|_{L^{p}}=\|u^{\ast}\|_{L^{p}}$.
\begin{definition}[\bf Weak $L^{p}$ spaces]\label{weakLp}
For $1\le p<\infty$, the weak $L^{p}$ (or Marcinkiewicz) space is defined as the space of all
measurable functions $u$ such that
 \begin{equation}
   \|u\|_{L^{p}_{w}}:=\sup_{\l>0}\l\m_{u}(\l)^{\frac{1}{p}}<\infty
 \end{equation}
 and $L^{\infty}_{w}:=L^{\infty}$.
\end{definition} Clearly $L^{p}\subseteq L^{p}_{w}$ and in
general the inclusion is strict: for instance, $u(x)=|x|^{-N/p}\in L^{p}_{w}(\R^{N})$ but does not
belong to any $L^{q}$. On the other hand 
\begin{equation}\label{ae24}
  L^{p}_{w}\cap L^{q}_{w}\subseteq L^{r},\qquad 1\le p<r<q\le\infty.
\end{equation}
\begin{definition}[\bf Lorentz spaces]\label{Lorentzs}
For $1\le p<\infty$ and $1\le q \le \infty$, the {Lorentz space} $L^{p,q}$ is defined as the set
of all measurable functions $u$ such that the following quasi-norm
\begin{align}
  \|u\|_{L^{p,q}}&:=\|t^{1/p}u^{\ast}(t)\|_{L^{q}_{\ast}}=
  \begin{cases}
    \left(\int\limits_{0}^{\infty}(t^{1/p}u^{\ast}(t))^q \frac{dt}{t}\right)^{\frac 1q} & \text{if }1\le q<\infty, \\
    \sup\limits_{t\ge0}t^{1/p}u^{\ast}(t) & \text{if }q=\infty,
  \end{cases}
\end{align}
is finite. We also set $L^{\infty,\infty}=L^{\infty}$.
\end{definition}
By H\"older's inequality
$L^{p,q_{1}}\subseteq L^{p,q_{2}}$ if $q_{1}\le q_{2}$ and more generally we have
\begin{equation}
  L^{p}\equiv L^{p,p}\subseteq L^{p,q}\subseteq L^{p,\infty}\equiv L^{p}_{w},\qquad 1\le p\le q\le \infty.
\end{equation}
Lorentz spaces have a classical characterization as interpolation of $L^{p}$ spaces (cf.
\cite{MR2424078}, Corollary 7.27).
\begin{proposition}\label{al2}
For any $1\le p_1<p<p_2\le \infty$ and $1\le q\le \infty$, we have
\begin{equation}\label{Lorentz_K2}
  L^{p,q}=\left(L^{p_1},L^{p_2}\right)_{\theta,q},\qquad
  \frac{1}{p}=(1-\theta)\frac{1}{p_1}+\theta\frac{1}{p_2},
\end{equation}
and in particular
\begin{equation}\label{Lorentz_K}
  L^{p,q}=\left(L^1,L^{\infty}\right)_{1-\frac 1p,q}.
\end{equation}
\end{proposition}
Another characterization of Lorentz spaces has been provided by L. Tartar with the aim of studying
improved Sobolev's embedding theorems.
\begin{lemma}[\cite{MR2328004}, Lemma 29.4]\label{lem_Lorentz}
Given a measurable function $u$ on $\R^N$ such that $\m_{u}(\l)<\infty$ for any $\l>0$, we
consider a {decreasing} sequence $(a_{k})_{k\in\Z}$ such that
\begin{equation}\label{ae33}
  u^{\ast}(e^{k})\le a_{k}\le u^{\ast}(e^{k}-),\qquad k\in\Z.
\end{equation}
Then, for ${1\le p}<\infty$ and $1\le q \le \infty$, we have
  $$u\in L^{p,q}\quad \text{ if and only if }\quad e^{k/p}a_{k} \in \ell^q(\Z).$$
Moreover, if $a_{k}\rightarrow 0$ as $k\rightarrow \infty$, then
  $$u\in L^{p,q} \quad \text{ if and only if }\quad e^{k/p}(a_{k}-a_{k+1}) \in \ell^q(\Z).$$
\end{lemma}
Notice that from \eqref{ae33} it follows that
\begin{equation}\label{ae27}
  \m_{u}(a_{k})\le e^{k}\le \m_{u}(a_{k}-)\le \m_{u}(a_{k+1}),\qquad k\in\Z.
\end{equation}

\bibliographystyle{acm}
\bibliography{bib1}

\def\cprime{$'$} \def\cprime{$'$} \def\cprime{$'$}
  \def\lfhook#1{\setbox0=\hbox{#1}{\ooalign{\hidewidth
  \lower1.5ex\hbox{'}\hidewidth\crcr\unhbox0}}} \def\cprime{$'$}
  \def\cprime{$'$} \def\cprime{$'$} \def\cprime{$'$} \def\cprime{$'$}
  \def\polhk#1{\setbox0=\hbox{#1}{\ooalign{\hidewidth
  \lower1.5ex\hbox{`}\hidewidth\crcr\unhbox0}}}
\begin{thebibliography}{10}

\bibitem{MR3951695}
{\sc Abedin, F., and Tralli, G.}
\newblock Harnack inequality for a class of {K}olmogorov-{F}okker-{P}lanck
  equations in non-divergence form.
\newblock {\em Arch. Ration. Mech. Anal. 233}, 2 (2019), 867--900.

\bibitem{MR2424078}
{\sc Adams, R.~A., and Fournier, J. J.~F.}
\newblock {\em Sobolev spaces}, second~ed., vol.~140 of {\em Pure and Applied
  Mathematics (Amsterdam)}.
\newblock Elsevier/Academic Press, Amsterdam, 2003.

\bibitem{MR1949176}
{\sc Bouchut, F.}
\newblock Hypoelliptic regularity in kinetic equations.
\newblock {\em J. Math. Pures Appl. (9) 81}, 11 (2002), 1135--1159.

\bibitem{MR2729292}
{\sc Bramanti, M., Cupini, G., Lanconelli, E., and Priola, E.}
\newblock Global {$L^p$} estimates for degenerate {O}rnstein-{U}hlenbeck
  operators.
\newblock {\em Math. Z. 266}, 4 (2010), 789--816.

\bibitem{MR4290567}
{\sc Camellini, F., Eleuteri, M., and Polidoro, S.}
\newblock A compactness result for the {S}obolev embedding via potential
  theory.
\newblock In {\em Harnack inequalities and nonlinear operators}, vol.~46 of
  {\em Springer INdAM Ser.} Springer, Cham, (2021), pp.~61--91.

\bibitem{MR3778645}
{\sc Cameron, S., Silvestre, L., and Snelson, S.}
\newblock Global a priori estimates for the inhomogeneous {L}andau equation
  with moderately soft potentials.
\newblock {\em Ann. Inst. H. Poincar\'{e} C Anal. Non Lin\'{e}aire 35}, 3
  (2018), 625--642.

\bibitem{MR2434049}
{\sc Di~Francesco, M., Pascucci, A., and Polidoro, S.}
\newblock The obstacle problem for a class of hypoelliptic ultraparabolic
  equations.
\newblock {\em Proc. R. Soc. Lond. Ser. A Math. Phys. Eng. Sci. 464}, 2089
  (2008), 155--176.

\bibitem{DiFrancescoPolidoro}
{\sc Di~Francesco, M., and Polidoro, S.}
\newblock Schauder estimates, {H}arnack inequality and {G}aussian lower bound
  for {K}olmogorov-type operators in non-divergence form.
\newblock {\em Adv. Differential Equations 11}, 11 (2006), 1261--1320.

\bibitem{MR4444079}
{\sc Dong, H., and Yastrzhembskiy, T.}
\newblock Global {$L_p$} estimates for kinetic {K}olmogorov-{F}okker-{P}lanck
  equations in nondivergence form.
\newblock {\em Arch. Ration. Mech. Anal. 245}, 1 (2022), 501--564.

\bibitem{MR657581}
{\sc Folland, G.~B., and Stein, E.~M.}
\newblock {\em Hardy spaces on homogeneous groups}, vol.~28 of {\em
  Mathematical Notes}.
\newblock Princeton University Press, Princeton, N.J.; University of Tokyo
  Press, Tokyo, 1982.

\bibitem{MR4444114}
{\sc Garofalo, N., and Tralli, G.}
\newblock Hardy-{L}ittlewood-{S}obolev inequalities for a class of
  non-symmetric and non-doubling hypoelliptic semigroups.
\newblock {\em Math. Ann. 383}, 1-2 (2022), 1--38.

\bibitem{MR3923847}
{\sc Golse, F., Imbert, C., Mouhot, C., and Vasseur, A.~F.}
\newblock Harnack inequality for kinetic {F}okker-{P}lanck equations with rough
  coefficients and application to the {L}andau equation.
\newblock {\em Ann. Sc. Norm. Super. Pisa Cl. Sci. (5) 19}, 1 (2019), 253--295.

\bibitem{lanpol}
{\sc Lanconelli, E., and Polidoro, S.}
\newblock On a class of hypoelliptic evolution operators.
\newblock {\em Rend. Sem. Mat. Univ. Politec. Torino 52}, 1 (1994), 29--63.

\bibitem{MR3726909}
{\sc Leoni, G.}
\newblock {\em A first course in {S}obolev spaces}, second~ed., vol.~181 of
  {\em Graduate Studies in Mathematics}.
\newblock American Mathematical Society, Providence, RI, 2017.

\bibitem{MR3753604}
{\sc Lunardi, A.}
\newblock {\em Interpolation theory}, vol.~16 of {\em Lecture Notes. Scuola
  Normale Superiore di Pisa}.
\newblock Edizioni della Normale, Pisa, 2018.

\bibitem{MR1751429}
{\sc Manfredini, M.}
\newblock The {D}irichlet problem for a class of ultraparabolic equations.
\newblock {\em Adv. Differential Equations 2}, 5 (1997), 831--866.

\bibitem{MR1762582}
{\sc Morbidelli, D.}
\newblock Fractional {S}obolev norms and structure of
  {C}arnot-{C}arath\'{e}odory balls for {H}\"{o}rmander vector fields.
\newblock {\em Studia Math. 139}, 3 (2000), 213--244.

\bibitem{MR3429628}
{\sc Pagliarani, S., Pascucci, A., and Pignotti, M.}
\newblock Intrinsic {T}aylor formula for {K}olmogorov-type homogeneous groups.
\newblock {\em J. Math. Anal. Appl. 435}, 2 (2016), 1054--1087.

\bibitem{Pascucci2}
{\sc Pascucci, A.}
\newblock H\"older regularity for a {K}olmogorov equation.
\newblock {\em Trans. Amer. Math. Soc. 355}, 3 (2003), 901--924.

\bibitem{MR2791231}
{\sc Pascucci, A.}
\newblock {\em P{DE} and martingale methods in option pricing}, vol.~2 of {\em
  Bocconi \& Springer Series}.
\newblock Springer, Milan; Bocconi University Press, Milan, 2011.

\bibitem{MR4355925}
{\sc Pascucci, A., and Pesce, A.}
\newblock On stochastic {L}angevin and {F}okker-{P}lanck equations: the
  two-dimensional case.
\newblock {\em J. Differential Equations 310\/} (2022), 443--483.

\bibitem{PascucciPolidoro}
{\sc Pascucci, A., and Polidoro, S.}
\newblock The {M}oser's iterative method for a class of ultraparabolic
  equations.
\newblock {\em Commun. Contemp. Math. 6}, 3 (2004), 395--417.

\bibitem{P22}
{\sc Pesce, A.}
\newblock {A}pproximation and {I}nterpolation in {K}olmogorov-type groups.
\newblock {\em arXiv:2205.05340. To appear in J. Math. Anal. Appl.\/} (2022).

\bibitem{MR1638177}
{\sc Polidoro, S., and Ragusa, M.~A.}
\newblock Sobolev-{M}orrey spaces related to an ultraparabolic equation.
\newblock {\em Manuscripta Math. 96}, 3 (1998), 371--392.

\bibitem{MR436223}
{\sc Rothschild, L.~P., and Stein, E.~M.}
\newblock Hypoelliptic differential operators and nilpotent groups.
\newblock {\em Acta Math. 137}, 3-4 (1976), 247--320.

\bibitem{MR558675}
{\sc Saka, K.}
\newblock Besov spaces and {S}obolev spaces on a nilpotent {L}ie group.
\newblock {\em Tohoku Math. J. (2) 31}, 4 (1979), 383--437.

\bibitem{MR0290095}
{\sc Stein, E.~M.}
\newblock {\em Singular integrals and differentiability properties of
  functions}.
\newblock Princeton Mathematical Series, No. 30. Princeton University Press,
  Princeton, N.J., 1970.

\bibitem{MR1662313}
{\sc Tartar, L.}
\newblock Imbedding theorems of {S}obolev spaces into {L}orentz spaces.
\newblock {\em Boll. Unione Mat. Ital. Sez. B Artic. Ric. Mat. (8) 1}, 3
  (1998), 479--500.

\bibitem{MR2328004}
{\sc Tartar, L.}
\newblock {\em An introduction to {S}obolev spaces and interpolation spaces},
  vol.~3 of {\em Lecture Notes of the Unione Matematica Italiana}.
\newblock Springer, Berlin; UMI, Bologna, 2007.

\bibitem{MR0216286}
{\sc Trudinger, N.~S.}
\newblock On imbeddings into {O}rlicz spaces and some applications.
\newblock {\em J. Math. Mech. 17\/} (1967), 473--483.

\bibitem{MR1218884}
{\sc Varopoulos, N.~T., Saloff-Coste, L., and Coulhon, T.}
\newblock {\em Analysis and geometry on groups}, vol.~100 of {\em Cambridge
  Tracts in Mathematics}.
\newblock Cambridge University Press, Cambridge, 1992.

\bibitem{ZhangBesov}
{\sc Zhang, X., and Zhang, X.}
\newblock Cauchy problem of stochastic kinetic equations.
\newblock {\em arXiv:2103.02267v1\/} (2022).

\end{thebibliography}

\end{document}